\documentclass[12pt, reqno]{amsart}
\usepackage{amsmath}
\usepackage{amssymb}
\usepackage{amsfonts}
\usepackage{amsthm}
\usepackage{amscd}
\usepackage{bm}
\usepackage{ulem}
\usepackage{graphicx}
\usepackage[dvipdfmx]{color}
\newtheorem{thrm}{Theorem}

\newtheorem{prob}{Problem}

\newtheorem{thm}{Theorem}[section]
\newtheorem{lem}[thm]{Lemma}
\newtheorem{pro}[thm]{Proposition}
\newtheorem{cor}[thm]{Corollary}

\newtheorem{conj}[thm]{Conjecture}
\newtheorem{defn}[thm]{Definition}

\newcommand{\Z}{\mathbf{Z}}
\newcommand{\Q}{\mathbf{Q}}
\newcommand{\N}{\mathbf{N}}

\begin{document}
\title[the basis-conjugating automorphism groups of free groups]
{On the structures of the Johnson cokernels of \\ the basis-conjugating automorphism groups of \\ free groups}
\author{Naoya Enomoto}
\address{Naoya Enomoto; The University of Electro-Communications,
1-5-1, Chofugaoka, Chofu city, Tokyo 182-8585, Japan.}
\email{enomoto-naoya@uec.ac.jp}
\author{Takao Satoh}
\address{Takao Satoh; Department of Mathematics, Faculty of Science Division II, Tokyo University of Science,
1-3, Kagurazaka, Shinjuku-ku, Tokyo 162-8601, Japan.}
\email{takao@rs.tus.ac.jp}
\subjclass[2020]{20F28(Primary), 20F14, 20C30, 20J06(Secondary)}
\keywords{Basis-conjugating automorphism groups of free groups, Johnson homomorphisms, Andreadakis-Johnson filtration.}
\maketitle
\begin{abstract}
In this paper, we study the Johnson homomorphisms 
of basis-conjugating automorphism groups of free groups.
We construct obstructions for the surjectivity of the Johnson homomorphisms.
By using it, we determine its cokernels of degree up to four, and give further observations for degree greater than four.

As applications, we give the affirmative answer for the Andreadakis problem for degree four,
and, we calculate the twisted first cohomology groups of the braid-permutation automorphism groups of a free group.

\end{abstract}
\tableofcontents
\section{Introduction}

Highly motivated by the study of the Johnson homomorphisms of the mapping class groups of surfaces
by pioneering works of Johnson \cite{Jo1}, Morita \cite{Mo1} and Hain \cite{Ha1} in the 1980s and 1990s, those of the automorphism groups of free groups
was initiated by several authors including Kawazumi \cite{Kaw} in the early 2000s.
Over the last two decades, good progress was
made in the study of the Johnson homomorphisms of the automorphism groups of free groups
through the works of many authors. 
In particular, the stable cokernel of the Johnson homomorphisms of the automorphism groups of free groups
have completely determined due to Satoh \cite{S03, S13} and Darn\'{e} \cite{Dar}. Enomoto-Satoh \cite{ES1} also gave a combinatorial description
of the GL-decomposition of the stable rational cokernels of the Johnson homomorphisms.
As is well-known due to a classical work of Artin, the automorphism group of a free group
contains the braid group of a plane. (See also Section 1.4 in \cite{Bir}.) 
Hence we have natural problems to determine the images and the cokernels of Johnson 
homomorphisms of the braid groups.
However, in contrast to the case of the mapping class group of a surface, there are few studies for them.
In general, it seems to be a difficult problem to attack this problem directly from a combinatorial group theretic and a representation theretic viewpoints.
In the present paper, we forcus on the basis-conjugating automorphism groups of free groups which contains the pure braid groups as a subgroup,
and study the Johsnon homomorphisms of them.

\vspace{0.5em}

Let us fix some notations. For any $n \geq 1$, let $F_n$ be the free group of rank $n$ with basis $x_1, \ldots , x_n$,
and $\mathrm{Aut}\,F_n$ the automorphism group of $F_n$. Set $H:=H_1(F_n,\Z)$ and $H^* := \mathrm{Hom}_{\Z}(H,\Z)$.
Denote by $\mathcal{L}_n=\bigoplus_{k \geq 1} \mathcal{L}_n(k)$ the free Lie algebra generated by $H$.
Let $\mathcal{A}_n(1) \supset \mathcal{A}_n(2) \supset \cdots$ be the Andreadakis-Johnson filtration of $\mathrm{Aut}\,F_n$.
Then for each $k \geq 1$, the $k$-th Johnson homomorphism
\[ \tau_k : \mathrm{gr}^k(\mathcal{A}_n):= \mathcal{A}_n(k)/\mathcal{A}_n(k+1) \rightarrow H^* \otimes_{\Z} \mathcal{L}_n(k+1) \]
is defined. Each $\tau_k$ is an injective $\mathrm{GL}(n,\Z)$-equivariant homomorphism. (See Section {\rmfamily \ref{S-John}} for details.)
Let 
\[ \mathcal{C}_n(k) := H^{\otimes k}/ \langle a_1 \otimes a_2 \cdots \otimes a_k - a_k \otimes a_1 \cdots \otimes a_{k-1} \,|\, a_i \in H \rangle \]
be the quotient module of $H^{\otimes k}$ by the permutation action of the cyclic group $C_k$ of degree $k$ on the components.
By combining results by Satoh \cite{S03, S13} and by Darn\'{e} \cite{Dar}, it is known that the cokernel of $\tau_k$ is isomorphic to
$\mathcal{C}_n(k)$ for $n \geq k+2$. In order to describe the cokernel of $\tau_k$, we used the trace map
\[ \mathrm{Tr}_k : H^* \otimes_{\Z} \mathcal{L}_n(k+1) \rightarrow H^* \otimes_{\Z} H^{\otimes (k+1)} \rightarrow H^{\otimes k} \rightarrow \mathcal{C}_n(k) \]
defined as the compositions of maps including a contraction. This is a generalization of the Morita trace defined in \cite{Mo1}.
(See Section {\rmfamily \ref{S-John}} for details.) Then we have a exact sequence
\[ 0 \rightarrow \mathrm{gr}^k(\mathcal{A}_n) \xrightarrow{\tau_k} H^* \otimes_{\Z} \mathcal{L}_n(k+1) \xrightarrow{\mathrm{Tr}_k} \mathcal{C}_n(k) \rightarrow 0 \]
for $n \geq k+2$. Furtheremore, Enomoto-Satoh \cite{ES1} gave a combinatorial description
of the GL-decomposition of $\mathcal{C}_n^{\Q}(k) := \mathcal{C}_n(k) \otimes_{\Z} \Q$.
In the paper, we consider the basis-conjugating automorphism group as an analogue of these works.

\vspace{0.5em}

An automorphism $\sigma$ of $F_n$ such that $x_i^{\sigma}$ is conjugate to $x_i$ for each $1 \leq i \leq n$ is called a basis-conjugating
automorphism of $F_n$. Let $\mathrm{P}\Sigma_n$ be the subgroup of $\mathrm{Aut}\,F_n$ consisting of all basis-conjugating automorphisms of $F_n$.
The group $\mathrm{P}\Sigma_n$ is called the basis-conjugating automorphism group of $F_n$.
Since McCool \cite{McC} gave the first finite presentation for $\mathrm{P}\Sigma_n$, it is also called the McCool group.
So far, $\mathrm{P}\Sigma_n$ has been studied by many authors from several viewpoints.

\vspace{0.5em}

For any $k \geq 1$, set $\mathcal{P}_n(k):= \mathrm{P}\Sigma_n \cap \mathcal{A}_n(k)$.
We have an injective $\mathfrak{S}_n$-equivariant homomorphism
\[ \tau_k^P : \mathrm{gr}^k(\mathcal{P}_n) := \mathcal{P}_n(k)/\mathcal{P}_n(k+1) \rightarrow \mathfrak{p}_n(k), \]
where $\mathfrak{p}_n(k)$ is a certain submodule of $H^* \otimes_{\Z} \mathcal{L}_n(k+1)$. (See Section {\rmfamily \ref{S-FLA-DA}} for details.)
From a veiwpoint of the study of the Johnson homomorphisms of $\mathrm{P}\Sigma_n$, we consider
the Lie algebra $\mathfrak{p}_n:=\bigoplus_{k \geq 1} \mathfrak{p}_n(k)$ in \cite{S18}. But
we should remark that $\mathfrak{p}_n$ was introduced as the Lie algebra of tangential derivations
by Alekseev and Torossian in \cite{AlT} for their study for the Kashiwara-Vergne problem.
The homomorphism $\tau_k^P$ is $\mathfrak{S}_n$-equivariant, and is called the $k$-th Johnson homomorphism of $\mathrm{P}\Sigma_n$. 
Here $\mathfrak{S}_n$ is the symmetric group of degree $n$, and is considered as the
subgroup of $\mathrm{GL}(n,\Z)$ consisting of all of permutation matrices.
(See also Section {\rmfamily \ref{S-FLA-DA}}.)
We should remark that in this paper the target of $\tau_k^P$ is always understood to be $\mathfrak{p}_n(k)$, not $H^* \otimes_{\Z} \mathcal{L}_n(k+1)$.
Hence the cokernel $\mathrm{Coker}(\tau_k^P)$ of $\tau_k^P$ means the quotient module of $\mathfrak{p}_n(k)$ by the image of $\tau_k^P$.
In \cite{S18}, we showed that $\tau_1^P$ is surjective and the abelianization of $\mathrm{P}\Sigma_n$.
Ibrahim \cite{Ibr} determined the structure of $\mathrm{gr}^2(\mathcal{P}_n)$ and $\mathrm{gr}^3(\mathcal{P}_3)$.
In particular, he obtained the ranks of them.

\vspace{0.5em}

Compared with the case of $\tau_k$, to describe the cokernel of $\tau_k^P$ for $k \geq 2$ is much more complicated.
Let $\mathrm{Tr}_k^P : \mathfrak{p}_n(k) \rightarrow \mathcal{C}_n(k)$ be the restriction of $\mathrm{Tr}_k$ to $\mathfrak{p}_n(k)$.
By the definition of $\mathcal{P}_n(k)$, it is easy to see that for any $k \geq 2$,
$\langle x_1^k, \ldots, x_n^k \rangle \cap \mathrm{Im}(\mathrm{Tr}_k^P) = \{ 0 \}$.
Let $\overline{\mathrm{Tr}}_k^P : \mathfrak{p}_n(k) \rightarrow \overline{\mathcal{C}}_n(k)$ be the composition of
$\mathrm{Tr}_k^P$ and
the natural projection $\mathcal{C}_n(k) \rightarrow \overline{\mathcal{C}}_n(k) := \mathcal{C}_n(k)/\langle x_1^k, \ldots, x_n^k \rangle$.
Then we have
\begin{thrm}[Propositions {\rmfamily \ref{P-1}} and {\rmfamily \ref{P-2}}, and Corollary {\rmfamily \ref{C-5}}]\label{I-1}
For any $n \geq 3$, we have an $\mathfrak{S}_n$-equivariant exact sequences
\[ 0 \rightarrow \mathrm{gr}^k(\mathcal{P}_n) \xrightarrow{\tau_k^P} \mathfrak{p}_n(k) \xrightarrow{\overline{\mathrm{Tr}}_k^P} \overline{\mathcal{C}}_n(k) \rightarrow 0 \]
for $2 \leq k \leq 3$, and
\[ 0 \rightarrow \mathrm{gr}^4(\mathcal{P}_n) \xrightarrow{\tau_4^P} \mathfrak{p}_n(4) \xrightarrow{\overline{\mathrm{Tr}}_4^P} \overline{\mathcal{C}}_n(4)
    \rightarrow Q \rightarrow 0, \]
where $Q$ is isomorphic to a free abelian group of rank $n(n-1)/2$ as an $\mathfrak{S}_n$-module.
\end{thrm}
\noindent
This shows that $\mathrm{Im}(\tau_k^P) = \mathrm{Ker}(\overline{\mathrm{Tr}}_k^P)$ for $k \leq 4$ and
the trace map $\overline{\mathrm{Tr}}_k^P$ is not surjective in general.
\begin{prob}\label{Prob-1}
In order to describe the image of $\tau_k^P$, there are the following basic problems. For any $k \geq 5$,
\begin{enumerate}
\item[(1)] Determine the image of $\overline{\mathrm{Tr}}_k^P$.
\item[(2)] Determine the gap at the inclusion $\mathrm{Im}(\tau_k^P) \subset \mathrm{Ker}(\overline{\mathrm{Tr}}_k^P)$.
\end{enumerate}
\end{prob}
\noindent
In general, it seems to be difficult to determine whether or not $\mathrm{Im}(\tau_k^P)=\mathrm{Ker}(\overline{\mathrm{Tr}}_k^P)$.
In this paper we show that a certain $\mathfrak{S}_n$-invariant subspace of $\mathrm{Ker}(\overline{\mathrm{Tr}}_k^P)$ is contained in
$\mathrm{Im}(\tau_k^P)$.
More precisely, for any composition $\alpha=(\alpha_1,\alpha_2, \ldots ,\alpha_n)$ of $k$,
denote by $\mathfrak{p}_n(k,\alpha)$ the subspace of $\mathfrak{p}_n(k)$ generated by
all $x_i^* \otimes [x_{i_1},x_{i_2}, \ldots ,x_{i_k}, x_i]$ such that $(i_1,i_2, \ldots, i_k) \in p(\alpha)$.
\begin{thrm}(=Theorem {\rmfamily \ref{T-0530}})
For any $n \geq 3$, $k \geq 1$ and any
composition $\alpha=(\alpha_1,\alpha_2, \ldots ,\alpha_n)$ of $k$ such that $\alpha_j=1$ for some $1 \leq j \leq n$,
we have $\mathrm{Ker}(\overline{\mathrm{Tr}}_k^P|_{\mathfrak{p}_n(k,\alpha)}) \subset \mathrm{Im}(\tau_k^P)$.
\end{thrm}
\noindent
To consider the above problems, we introduce a variant $\widetilde{\mathrm{Tr}}_k^P$ of the trace map $\overline{\mathrm{Tr}}_k^P$.
For any $k \geq 2$,
let $N(k)$ be the $\mathfrak{S}_n$-invariant submodule of $\mathcal{C}_n(k)$ generated by
$x_{j_1} x_{j_2} \cdots x_{j_k}$ for all $1 \leq j_1, \ldots, j_k \leq n$ such that each $j_l$ $(1 \leq l \leq k)$ appears at least twice in
$j_1, j_2, \ldots, j_k$. Let $\widetilde{\mathrm{Tr}}_k^P : \mathfrak{p}_n(k) \rightarrow \widetilde{\mathcal{C}}_n(k):= \mathcal{C}_n(k)/N(k)$
be the composition of $\mathrm{Tr}_k^P$ and the natural projection $\mathcal{C}_n(k) \rightarrow \widetilde{\mathcal{C}}_n(k)$.
Then we have
\[ \mathrm{Im}(\tau_k^P) \subset \mathrm{Ker}(\mathrm{Tr}_k^P) = \mathrm{Ker}(\overline{\mathrm{Tr}}_k^P) \subset
   \mathrm{Ker}(\widetilde{\mathrm{Tr}}_k^P) \subset \mathfrak{p}_n(k). \]
\begin{thrm}[Theorem {\rmfamily \ref{T-main_1}}]\label{I-2}
For any $n \geq 2$ and $k \geq 2$, the map $\widetilde{\mathrm{Tr}}_k^P$ is surjective. 
\end{thrm}
\noindent
We give some observations for Problem {\rmfamily \ref{Prob-1}} in Section {\rmfamily \ref{S-OBS}}.
First, we find the explicite dimensions of $\mathrm{Im}(\overline{\mathrm{Tr}}_k^P)$ and 
$\mathrm{Ker}(\overline{\mathrm{Tr}}_k^P)$ for $5 \leq k \leq 9$ and general $n$. Second, we also give the explicite dimensions of the Johnson images 
$\mathrm{Im}(\tau_k^P)$ for $5 \leq k \leq 8$ and $n=3$.
We will see some gaps at the inclusion $\mathrm{Im}(\tau_k^P) \subsetneq \mathrm{Ker}(\overline{\mathrm{Tr}}_k^P)$ for $k=7$ and $8$, and $n=3$.

\vspace{0.5em}

The study of the Johnson homomorphisms gives us several applications. Here we consider three applications.
The first one is to apply to the Andreadakis problem for $\mathrm{P}\Sigma_n$.
Let $\mathrm{P}\Sigma_n=\mathrm{P}\Sigma_n{(1)} \supset \mathrm{P}\Sigma_n{(2)} \supset \cdots$ be the lower central series of $\mathrm{P}\Sigma_n$.
We have $\mathrm{P}\Sigma_n{(k)} \subset \mathcal{P}_n(k)$ for any $k \geq 1$.
The Andreadakis problem for $\mathrm{P}\Sigma_n$ asks whether $\mathrm{P}\Sigma_n{(k)} = \mathcal{P}_n(k)$ or not.
In general, except for some low degree or unstable cases, the Andreadakis problems for $\mathrm{P}\Sigma_n$ and $\mathrm{Aut}\,F_n$ are still open.
Satoh \cite{S15} showed that $\mathrm{P}\Sigma_n{(k)} = \mathcal{P}_n(k)$ for $k \leq 3$.
By using the thrid Johnson homomorphism, we obtain the following.
\begin{thrm}[Corollary {\rmfamily \ref{C-6}}]\label{I-3}
For $n \geq 3$, we have $\mathcal{P}_n(4)=\mathrm{P}\Sigma_n(4)$.
\end{thrm}
\noindent
This shows that Ibrahim's conjecture (See Conjecture A.1. in \cite{Ibr}) holds for any $n \geq 3$ and $k=3$.

\vspace{0.5em}

We remark that by using the results above,
we can reobtain results for the second homology and cohomology groups of $\mathrm{P}\Sigma_n$ obtained by
Brownstein and Lee \cite{BrL}, and Jensen, McCammond and Meiyer \cite{JMM}.
The first and the second integral homology groups of $\mathrm{P}\Sigma_n$ are determined by Brownstein and Lee \cite{BrL}, and
the integral cohomology ring of $\mathrm{P}\Sigma_n$ is determined by Jensen, McCammond and Meiyer \cite{JMM}.
By using the second Johnson homomorphism and the McCool presentation for $\mathrm{P}\Sigma_n$,
in Propositions \ref{T-main-4} and \ref{T-main-5} we obtain the following. For $n \geq 3$,
\begin{itemize}
\item[$(1)$]\,\,$H_2(\mathrm{P}\Sigma_n,\Z) \cong \Z^{\frac{1}{2}n^2(n-1)(n-2)}$. 
\item[$(2)$]\,\,The cup product map $\cup : \wedge^2 H^1(\mathrm{P}\Sigma_n, \Z) \rightarrow H^2(\mathrm{P}\Sigma_n, \Z)$ is surjective.
\end{itemize}

\vspace{0.5em}

Finally, we consider twisted first homology groups of the braid-permutation automorphism groups of free groups.
An element $\sigma \in \mathrm{Aut}\,F_n$ is a braid-permutation automorphism of $F_n$ if there exists some $\mu \in \mathfrak{S}_n$ and some
$a_i \in F_n$ for any $1 \leq i \leq n$ such that
\[ x_i^{\sigma}=a_i^{-1} x_{\mu(i)} a_i \hspace{1em} (1 \leq i \leq n). \]
Let $\mathrm{BP}_n$ be the subgroup of $\mathrm{Aut}\,F_n$ consisting of all braid-permutation automorphisms of $F_n$.
According to \cite{FRR}, we call $\mathrm{BP}_n$ the braid-permutation automorphism group of $F_n$.
We remark that $\mathrm{BP}_n$ is isomorphic to the loop braid group on $n$ components.
(For details, see Sections 3 and 4 in \cite{Dam} for example.)
We see $\mathrm{BP}_n \cap \mathcal{A}_n(1) = \mathrm{P}\Sigma_n$.
Let $V$ be the rational standard irreducible representation of $\mathfrak{S}_n$.
By observing the decomposition of the image of $\tau_{1, \Q}^P$,
we see that $V$ appears in the decomposition of $H_1(\mathrm{P}\Sigma_n, \Q)$
By using a finite presentation for $\mathrm{BP}_n$ obtained by Fenn-Rim\'{a}nyi-Rourke \cite{FRR},
as an analogue of calculations of twisted first cohomology groups of $\mathrm{Aut}\,F_n$ in \cite{S01}, we obtain the following.
\begin{thrm}[Theorem {\rmfamily \ref{T-main-6}}]\label{I-4}
For $n \geq 3$, $H^1(\mathrm{BP}_n, V) \cong \Z^{\oplus 2} \oplus \Z/4\Z$.
\end{thrm}
\noindent
Let $B_n$ be Artin's braid group of $n$ strands. Due to a classical work of Artin, it is known that $B_n$ can be regard as a subgroup of $\mathrm{Aut}\,F_n$.
From this viewpoint, by observing the calculations in the proof of Theorem {\rmfamily \ref{I-4}}, we can see that for $n \geq 3$,
\[ H^1(B_n,V) = \Z \oplus \Z/4\Z, \hspace{1em} H^1(\mathfrak{S}_n,V) = \Z/4\Z. \]

\tableofcontents

\section{Notation and conventions}\label{S-Not}

In this section, we fix notation and conventions.
In this paper, a group means a multiplicative group if otherwise noted.
Let $G$ be a group. 

\begin{itemize}
\item The abelianization of $G$ is denoted by $G^{\mathrm{ab}}$ if otherwise noted.
\item The automorphism group $\mathrm{Aut}\,G$ acts on $G$ from the right. For any $\sigma \in \mathrm{Aut}\,G$ and $x \in G$,
      the action of $\sigma$ on $x$ is denoted by $x^{\sigma}$.
\item For a normal subgroup $N$, we often denote the coset class of an element $g \in G$ by the same $g$ in the quotient group $G/N$ if
      there is no confusion.
\item For elements $x$ and $y$ of $G$, the commutator bracket $[x,y]$ of $x$ and $y$
      is defined to be $[x,y]:=xyx^{-1}y^{-1}$. Then for any $x, y, z \in G$, we have
\begin{equation}
  [xy,z] =[x,[y,z]][y,z][x,z], \hspace{1em} [x,yz] = [x,y][x,z][[z,x],y]  \label{eq-commutator_form_1}
\end{equation}
and
\begin{equation}
  [x^{-1},z]=[[x^{-1}, z],x][x,z]^{-1}, \hspace{1em} [x,y^{-1}] = [x,y]^{-1} [y,[y^{-1}, x]]. \label{eq-commutator_form_2}
\end{equation}
\item For subgroups $H$ and $K$ of $G$, we denote by $[H,K]$ the commutator subgroup of $G$
      generated by $[h, k]$ for all $h \in H$ and $k \in K$.
\item For elements $g_1, \ldots, g_k \in G$, a left-normed commutator
\[ [[ \cdots [[ g_{1},g_{2}],g_{3}], \cdots ], g_{k}] \in G \]
of weight $k$ is denoted by $[g_{1},g_{2}, \cdots, g_{k}]$.
\item For any $\Z$-module $A$, we denote the $\Q$-vector space $A \otimes_{\Z} \Q$ by the symbol obtained by attaching a subscript $\Q$ to $A$,
like $A_{\Q}$ or $A^{\Q}$. Similarly, for any $\Z$-linear map $f: A \rightarrow B$,
the induced $\Q$-linear map $A_{\Q} \rightarrow B_{\Q}$ is denoted by $f_{\Q}$ or $f^{\Q}$.
\end{itemize}

\section{Free Lie algebras and its derivation algebras}\label{S-FLA-DA}

In this section, we review a free Lie algebra and its derivation algebra.
Furthermore, we consider a certain Lie subalgebra of the derivation algebra of the free Lie algebra.

\vspace{0.5em}

For $n \geq 1$, let $F_n$ be the free group of rank $n$ with basis $x_1, \ldots , x_n$.
Let $\Gamma_n(1) \supset \Gamma_n(2) \supset \cdots$ be the lower central series of $F_n$ defined by
\[ \Gamma_n(1):= F_n, \hspace{1em} \Gamma_n(k) := [\Gamma_n(k-1),F_n], \hspace{1em} k \geq 2. \]
We denote by $\mathcal{L}_n(k) := \Gamma_n(k)/\Gamma_n(k+1)$ the $k$-th graded quotient of the lower central series of $F_n$,
and by $\mathcal{L}_n := {\bigoplus}_{k \geq 1} \mathcal{L}_n(k)$ the associated graded sum.
For each $k \geq 1$, we consider $\mathcal{L}_n(k)$ as an additive group. For example, for any $x, y, z \in F_n$
we have
\[ [xy,z]=[x,z]+[y,z], \hspace{1em} [x^{-1},y]=-[x,y] \]
in $\mathcal{L}_n(1)$ from (\ref{eq-commutator_form_1}) and (\ref{eq-commutator_form_2}).
The graded sum $\mathcal{L}_n$ has the Lie algebra structure with the Lie bracket induced from the commutator bracket of $F_n$.
It is known that the Lie algebra $\mathcal{L}_n$ is isomorphic to the free Lie algebra generated by $\mathcal{L}_n(1)$
by a classical work of Magnus. (For details, see \cite{MKS} and \cite{Reu} for example.)
Thus each of $\mathcal{L}_n(k)$ is a free abelian group of finite rank.
Witt \cite{Wit} gave a formula for the rank of $\mathcal{L}_n(k)$ as
\begin{equation}\label{ex-witt}
 r_n(k) := \frac{1}{k} \sum_{d | k} \mu(d) n^{\frac{k}{d}}
\end{equation}
where $\mu$ is the M\"{o}bius function, and $d$ runs over all positive divisors of $k$.
Hall \cite{Hl1} constructed an explicit basis of $\mathcal{L}_n(k)$ with basic commutators of weight $k$.
For example, basic commutators of weight less than five are listed below.
\vspace{0.5em}
\begin{center}
{\renewcommand{\arraystretch}{1.3}
\begin{tabular}{|c|l|l|} \hline
  $k$  & basic commutators  &     \\ \hline
  $1$  & $x_i$ & $1 \leq i \leq n$                 \\ \hline
  $2$  & $[x_i,x_j]$        & $1 \leq j<i \leq n$  \\ \hline
  $3$  & $[x_i, x_j, x_l]$  & $i > j \leq l$    \\ \hline
  $4$  & $[x_i, x_j, x_l, x_m]$  & $i > j \leq l \leq m$  \\
       & $[[x_i, x_j], [x_l, x_m]]$  & $(i,j) > (l,m)$ \\ \hline
\end{tabular}}
\end{center}
\vspace{0.5em}
\noindent
Here $(i,j) > (l,m)$ means the lexicographic order on $\N \times \N$.
(For details, see \cite{Hl2} and \cite{Reu} for example.)

\vspace{0.5em}

For simplicity, we denote by $H$ the abelianization $F_n^{\mathrm{ab}}=H_1(F_n,\Z)$ of $F_n$.
The basis $x_1, \ldots, x_n$ of $F_n$ induces a basis of $H$ as a free abelian group. we also denote it by the same letters $x_1, \ldots, x_n$
by abuse of language.
By fixing this basis of $H$, we identify the automorphism group $\mathrm{Aut}\,H$ with the general linear group $\mathrm{GL}(n,\Z)$.
Then $\mathrm{GL}(n,\Z)$ naturally acts on each of $\mathcal{L}_n(k)$.
For example, as a $\mathrm{GL}(n,\Z)$-module we have $\mathcal{L}_n(1) = H$ and $\mathcal{L}_n(2) = \wedge^2 H$.
Furthermore, $\mathrm{GL}(n,\Q)$ naturally acts on each of $\mathcal{L}_n^{\Q}(k)$.
The irreducible decompositions of $\mathcal{L}_n^{\Q}(k)$ for $1 \leq k \leq 4$ are given as follows.
\vspace{0.5em}
\begin{center}
{\renewcommand{\arraystretch}{1.3}
\begin{tabular}{|c|l|l|l|} \hline
  $k$  & $\mathcal{L}_n(k)$ & $\mathrm{rank}_{\Z}(\mathcal{L}_n(k))$  \\ \hline
  $1$  & $[1]$     & $n$   \\ \hline
  $2$  & $[1^2]$   & $n(n-1)/2$ \\ \hline
  $3$  & $[2,1]$   & $n(n^2-1)/3$ \\ \hline
  $4$  & $[3,1] \oplus [2,1^2]$ & $n^2(n^2-1)/4$ \\ \hline
\end{tabular}}
\end{center}
\vspace{0.5em}
\noindent
Here $[\lambda]$ means the irreducible polynomial representation of $\mathrm{GL}(n,\Q)$ associated to the Young tableau $\lambda$.

\vspace{0.5em}

The universal enveloping algebra of $\mathcal{L}_n$ is the tensor algebra.
\[ T(H):= \Z \oplus H \oplus H^{\otimes 2} \oplus \cdots \]
of $H$ over $\Z$. From Poincar\'{e}-Birkhoff-Witt's theorem, the natural map $\iota : \mathcal{L}_n \rightarrow T(H)$
is injective. We denote by $\iota_k : \mathcal{L}_n(k) \rightarrow H^{\otimes k}$ the homomorphism of degree $k$ part of $\iota$, and
consider $\mathcal{L}_n(k)$ as a $\mathrm{GL}(n,\Z)$-submodule $H^{\otimes k}$ through $\iota_k$.
Similarly, we consider $\mathcal{L}_n^{\Q}(k)$ as a $\mathrm{GL}(n,\Q)$-submodule $H_{\Q}^{\otimes k}$ through $\iota_k^{\Q}$.

\vspace{0.5em}

Next, we consider the derivation algebra of the free Lie algebra.
We denote the Lie algebra of derivations of $\mathcal{L}_n$ by
\[ \mathrm{Der}(\mathcal{L}_n) := \{ f : \mathcal{L}_n \xrightarrow{\Z-\mathrm{linear}} \mathcal{L}_n \,|\, f([x,y]) = [f(x),y]+ [x,f(y)],
    \,\,\, x, y \in \mathcal{L}_n \}. \]
The Lie algebra $\mathrm{Der}(\mathcal{L}_n)$ is a graded Lie algebra.
For $k \geq 0$, the degree $k$ part of $\mathrm{Der}(\mathcal{L}_n)$ is defined to be
\[ \mathrm{Der}(\mathcal{L}_n)(k) := \{ f \in \mathrm{Der}(\mathcal{L}_n) \,|\, f(x) \in \mathcal{L}_n(k+1), \,\,\, x \in H \}. \]
For $k \geq 1$, by considering the universality of the free Lie algebra, we indetify $\mathrm{Der}(\mathcal{L}_n)(k)$ with
$\mathrm{Hom}_{\Z}(H,\mathcal{L}_n(k+1)) = H^* {\otimes}_{\Z} \mathcal{L}_n(k+1)$
as a $\mathrm{GL}(n,\Z)$-module.
Let $x_1^*, \ldots, x_n^*$ be the dual basis of $x_1, \ldots, x_n \in H$. Then $\mathrm{Der}(\mathcal{L}_n)(k)$
is generated by $x_i^* \otimes [x_{j_1}, x_{j_2}, \ldots, x_{j_{k+1}}]$ for all $1 \leq i, j_1, \ldots, j_{k+1} \leq n$
as a $\Z$-module.
Let $\mathrm{Der}^+(\mathcal{L}_n)$ be the graded Lie subalgebra of $\mathrm{Der}(\mathcal{L}_n)(k)$ with positive degree.
Similarly, we define the graded Lie algebra $\mathrm{Der}^+(\mathcal{L}_{n}^{\Q})$.
Then, we have $\mathrm{Der}^+(\mathcal{L}_{n}^{\Q}) = \mathrm{Der}^+(\mathcal{L}_n) \otimes_{\Z} \Q$, and
$\mathrm{Der}(\mathcal{L}_n^{\Q})(k)=H_{\Q}^* {\otimes}_{\Q} \mathcal{L}_n^{\Q}(k+1)$ for any $k \geq 1$.
\begin{defn}
For any $k \geq 1$, let $\mathfrak{p}_n(k)$ be the $\Z$-submodule of $\mathrm{Der}(\mathcal{L}_n)(k)$ generated by elements
\[ x_i^* \otimes [x_{j_1}, \ldots, x_{j_k}, x_i] \]
for all $1 \leq i, j_1, \ldots, j_k \leq n$. Let $\mathfrak{p}_n$ be the graded sum $\bigoplus_{k \geq 1} \mathfrak{p}_n(k)$.
\end{defn}
\noindent
Alekseev and Torossian \cite{AlT} showed that $\mathfrak{p}_n$ is a Lie subalgebra of $\mathrm{Der}^+(\mathcal{L}_n)$.
(See also Satoh \cite{S18}.)
We remark that
$x_i^* \otimes [x_j, x_i]$ for any $1 \leq i \neq j \leq n$ form a basis of $\mathfrak{p}_n(1)$, and
$\mathfrak{p}_n(k)$ is isomorphic to $H \otimes_{\Z} \mathcal{L}_n(k)$ as an $\mathfrak{S}_n$-module.
Hence we have $\mathrm{rank}_{\Z}(\mathfrak{p}_n(k)) = nr_n(k)$ for any $k \geq 2$.

\section{Johonson homomorphisms}\label{S-John}

Here we recall the Andreadakis-Johnson filtration and the Johnson homomorphisms of the automorphism groups of free groups.
For each $k \geq 1$, the action of $\mathrm{Aut}\,F_n$ on the nilpotent quotient group $F_n/\Gamma_n(k+1)$ of $F_n$ induces the homomorphism
\[ \mathrm{Aut}\,F_n \rightarrow \mathrm{Aut}(F_n/\Gamma_n(k+1)). \]
We denote its kernel by $\mathcal{A}_n(k)$. Then the groups $\mathcal{A}_n(k)$ define the descending filtration
\[ \mathcal{A}_n(1) \supset \mathcal{A}_n(2) \supset \mathcal{A}_n(3) \supset \cdots \]
of $\mathrm{Aut}\,F_n$.
This filtration is called the Andreadakis-Johnson filtration of $\mathrm{Aut}\,F_n$, and the first term is called the IA-automorphism group
of $F_n$, denoted by $\mathrm{IA}_n$.
Historically, the Andreadakis-Johnson filtration was originally introduced by Andreadakis \cite{And} in the 1960s.
He showed that
\begin{thm}[Andreadakis \cite{And}]\label{T-And} \quad
\begin{enumerate}
\item For any $k$, $l \geq 1$, $\sigma \in \mathcal{A}_n(k)$ and $x \in \Gamma_n(l)$, $x^{-1} x^{\sigma} \in \Gamma_n(k+l)$.
\item For any $k$ and $l \geq 1$, $[\mathcal{A}_n(k), \mathcal{A}_n(l)] \subset \mathcal{A}_n(k+l)$.
\item $\displaystyle \bigcap_{k \geq 1} \mathcal{A}_n(k) =\{ 1 \}$.
\end{enumerate}
\end{thm}
\noindent
For each $k \geq 1$, the group $\mathrm{Aut}\,F_n$ acts on $\mathcal{A}_n(k)$ by conjugation, and
it naturally induces the action of $\mathrm{GL}(n,\Z)=\mathrm{Aut}\,F_n/\mathrm{IA}_n$ on
the graded quotients $\mathrm{gr}^k (\mathcal{A}_n) := \mathcal{A}_n(k)/\mathcal{A}_n(k+1)$ by Part (2) of Theorem {\rmfamily \ref{T-And}}.

\vspace{0.5em}

In order to study the $\mathrm{GL}(n,\Z)$-module structure of ${\mathrm{gr}}^k (\mathcal{A}_n)$,
we consider the Johnson homomorphisms.
For each $k \geq 1$, define the homomorphism
$\tilde{\tau}_k : \mathcal{A}_n(k) \rightarrow \mathrm{Hom}_{\Z}(H, {\mathcal{L}}_n(k+1))$ by
\[ \sigma \hspace{0.3em} \mapsto \hspace{0.3em} (x \mapsto x^{-1} x^{\sigma}), \hspace{1em} x \in H. \]
The kernel of $\tilde{\tau}_k$ is $\mathcal{A}_n(k+1)$. 
Hence, it induces the injective homomorphism
\[ \tau_k : \mathrm{gr}^k (\mathcal{A}_n) \hookrightarrow \mathrm{Hom}_{\Z}(H, \mathcal{L}_n(k+1))
       = H^* \otimes_{\Z} \mathcal{L}_n(k+1). \]
The homomorphisms $\tilde{\tau}_k$ and ${\tau}_{k}$ are called the $k$-th Johnson homomorphisms of $\mathrm{Aut}\,F_n$.
Each of ${\tau}_{k}$ is $\mathrm{GL}(n,\Z)$-equivariant.
The graded sum $\mathrm{gr}(\mathcal{A}_n) := \bigoplus_{k \geq 1} \mathrm{gr}^k (\mathcal{A}_n)$ has the Lie algebra structure
with the Lie bracket induced from the commutator bracket on $\mathrm{IA}_n$.
Furthermore, the graded sum $\tau := \bigoplus_{k \geq 1} \tau_k : \mathrm{gr}(\mathcal{A}_n) \rightarrow \mathrm{Der}^+(\mathcal{L}_n)$
is a Lie algebra homomorphism.

\vspace{0.5em}

Here we remark on the cokernels of the Johnson homomorphisms.
For $k \geq 1$, let ${\varphi}^{k} : H^* {\otimes}_{\Z} H^{\otimes (k+1)} \rightarrow H^{\otimes k}$
be the contraction map defined by
\[ x_i^* \otimes x_{j_1} \otimes \cdots \otimes x_{j_{k+1}} \mapsto x_i^*(x_{j_1}) \, \cdot
    x_{j_2} \otimes \cdots \otimes \cdots \otimes x_{j_{k+1}}. \]
For the natural embedding ${\iota}_{k+1} : \mathcal{L}_n(k+1) \rightarrow H^{\otimes (k+1)}$,
we obtain the $\mathrm{GL}(n,\Z)$-equivariant homomorphism
\[ \Phi^k = {\varphi}^{k} \circ ({id}_{H^*} \otimes {\iota}_{k+1})
    : H^* {\otimes}_{\Z} \mathcal{L}_n(k+1) \rightarrow H^{\otimes k}. \]
We also call $\Phi^k$ the contraction map of degree $k$.
From Lemma 3.1 in \cite{S03}, we have the following.
\begin{lem}\label{L-S03_1}
For any $1 \leq i, i_1, \ldots, i_{k} \leq n$ such that $i_1 \neq i$, not required $i_1, \ldots, i_k$ are distinct, we have
\[\begin{split}
   \Phi^k(x_i^* \otimes & [x_i, x_{i_1}, \ldots, x_{i_{k}}]) \\
     & = x_{i_1} \otimes \cdots \otimes x_{i_{k}}
         - \sum_{l=2}^k \delta_{i, i_l} [x_i, x_{i_1}, \ldots, x_{i_{l-1}}] \otimes x_{i_{l+1}} \otimes \cdots \otimes x_{i_k}.
  \end{split} \]
\end{lem}

\vspace{0.5em}

Let $\mathcal{C}_n(k)$ be the quotient module of $H^{\otimes k}$ by the action of the cyclic group $C_k$ of order $k$ on the components.
Namely,
\[ \mathcal{C}_n(k) := H^{\otimes k} \big{/} \langle a_1 \otimes a_2 \otimes \cdots \otimes a_k - a_2 \otimes a_3 \otimes \cdots \otimes a_k \otimes a_1
   \,|\, a_i \in H \rangle. \]
Let $\varpi^k : H^{\otimes k} \rightarrow \mathcal{C}_n(k)$ be the natural projection.
For any $1 \leq j_1, \ldots, j_k \leq n$, we write the image of $x_{j_1} \otimes \cdots \otimes x_{j_k} \in H^{\otimes k}$
by the map $\varpi^k$ as $x_{j_1} x_{j_2} \cdots x_{j_k} \in \mathcal{C}_n(k)$. Thus, we have
$x_{j_1} x_{j_2} \cdots x_{j_k}=x_{j_2} \cdots x_{j_k} x_{j_1}$ in $\mathcal{C}_n(k)$.
By \cite{Reu}, it is known that $\mathcal{C}_n(k)$ is a free abelian group, and
the necklaces among $x_1, \ldots, x_n$ of length $k$ form a basis of $\mathcal{C}_n(k)$.
Hence we have
\[ \mathrm{rank}_{\Z}(\mathcal{C}_n(k)) = \frac{1}{k}\sum_{d|k}\varphi(d)n^{\frac{k}{d}} \]
where $\varphi : \N \rightarrow \N$ is the Euler function, and $d$ runs over all positive divisors of $k$.

\vspace{0.5em}

Set
$\mathrm{Tr}_k := \varpi^k \circ \Phi^k : H^* {\otimes}_{\Z} \mathcal{L}_n(k+1) \rightarrow \mathcal{C}_n(k)$.
We call $\mathrm{Tr}^k$ the trace map of degree $k$.
By combining the results by Satoh \cite{S03, S13} and Darne \cite{Dar}, we see that for any $k \geq 2$ and $n \geq k+2$, the sequence
\begin{equation}
 0 \rightarrow \mathrm{gr}^k (\mathcal{A}_n) \xrightarrow{\tau_k} H^* {\otimes}_{\Z} \mathcal{L}_n(k+1) \xrightarrow{\mathrm{Tr}_k} \mathcal{C}_n(k) \rightarrow 0 \label{eq-John-trace}
\end{equation}
is a $\mathrm{GL}(n,\Z)$-equivariant exact sequence. Furthermore, Enomoto-Satoh \cite{ES1} gave a combinatorial description of the $\mathrm{GL}(n,\Q)$-irreducible decomposition
of the stable rational cokernel $\mathcal{C}_n^{\Q}(k)$ of the Johnson homomorphism $\tau_{k,\Q}$. For $1 \leq k \leq 5$, the decompositions are given as follows.
\vspace{0.5em}
\begin{center}
{\renewcommand{\arraystretch}{1.3}
\begin{tabular}{|c|l|l|} \hline
  $k$  & $\mathcal{C}_n^\Q(k)$ & $\mathrm{dim}_{\Q}(\mathcal{C}_n^\Q(k))$ \\ \hline
  $2$  & $[2]$                 & $n(n+1)/2$ \\ \hline
  $3$  & $[3] \oplus [1^3]$    & $n(n^2+2)/3$  \\ \hline
  $4$  & $[4] \oplus [2,2] \oplus [2,1^2]$ & $n(n+1)(n^2-n+2)/4$ \\ \hline
  $5$  & $[5] \oplus [3,2] \oplus 2[3,1^2] \oplus [2^2,1] \oplus [1^5]$ & $n(n^4+4)/5$  \\ \hline
\end{tabular}}
\end{center}
\vspace{0.5em}
\noindent
Here we remark that for any $k \geq 2$, the irreducible component $[k]$, which is the symmetric product of $H_{\Q}$ of degree $k$,
appears in $\mathcal{C}_n^{\Q}(k)$ with multiplicity one. This result was originally obtained by Morita \cite{Mo1}.
Enomoto-Satoh \cite{ES1} showed that $[1^k]$ for any odd $k \geq 3$ and $n \geq k+2$, and $[2,1^{k-2}]$ for even $k \geq 4$ and $n \geq k+2$
appear in $\mathcal{C}_n^{\Q}(k)$ with multiplicity one. 

\vspace{0.5em}

Here we remark the Andreadakis problem on $\mathrm{Aut}\,F_n$.
Let $\mathrm{IA}_n = \mathcal{A}_n'(1) \supset \mathcal{A}_n'(2) \supset \cdots$ be the lower central series of $\mathrm{IA}_n$,
Since the Andreadakis-Johnson filtration is a central filtration of $\mathrm{IA}_n$,
we see $\mathcal{A}_n'(k) \subset \mathcal{A}_n(k)$ for any $k \geq 1$.
\begin{conj}[The Andreadakis conjecture]
For any $n \geq 3$ and $k \geq 1$, $\mathcal{A}_n'(k) = \mathcal{A}_n(k)$.
\end{conj}
\noindent
Andreadakis showed that the above conjecture is true for $n=2$ and any $k \geq 2$, and $n=3$ and $k \leq 3$.
It is known that $\mathcal{A}_n'(2) = \mathcal{A}_n(2)$ for any $n \geq 2$ due to Bachmuth \cite{Ba2}.
This result is also obtained from the fact that the first Johnson homomorphism is the abelianization of $\mathrm{IA}_n$ by independent works
of Cohen-Pakianathan \cite{Co1, Co2}, Farb \cite{Far} and Kawazumi \cite{Kaw}.
Satoh \cite{S25} showed that $\mathcal{A}_n'(3) = \mathcal{A}_n(3)$ for $n \geq 3$.
Bartholdi \cite{LB1, LB2} showed that this conjecture is not true in general. In particular, he showed that
$\mathcal{A}_3(4)/\mathcal{A}_3'(4) \cong (\Z/2\Z)^{\oplus 14} \oplus (\Z/3\Z)^{\oplus 9}$ and
$\mathcal{A}_3(5)/\mathcal{A}_3'(5) \otimes_{\Z} \Q \cong \Q^{\oplus 3}$.
In the stable range, the Andreadakis conjecture is still open problem.
We remark that Darn\'{e} \cite{Dar} gave affirmative answers for the Andreadakis problems restricted to certain subgroups of $\mathrm{Aut}\,F_n$.

\section{The basis-conjugating automorphism groups}

Here we recall the basis-conjugating automorphisms of $F_n$.
Let $\mathrm{P}\Sigma_n$ be the subgroup of $\mathrm{Aut}\,F_n$ consisting of all basis-conjugating automorphisms of $F_n$.
The group $\mathrm{P}\Sigma_n$ is a subgroup of $\mathrm{IA}_n$. For any $1 \leq i \neq j \leq n$, let $K_{ij}$ be the basis-conjugating automorphism of $F_n$ defined by
\[ x_t \mapsto \begin{cases}
               {x_j}^{-1} x_i x_j & t=i, \\
               x_t                & t \neq i.
              \end{cases}\]
McCool gave a finite presentation for $\mathrm{P}\Sigma_n$ as follows.
\begin{thm}[McCool \cite{McC}]
The group $\mathrm{P}\Sigma_n$ is generated by $K_{ij}$ for $1 \leq i \neq j \leq n$ subject to relations: \\
\hspace{2em} {\bf{(P1)}}: $[K_{ij},K_{kj}]=1$ for $(i,j) < (k,j)$, \\
\hspace{2em} {\bf{(P2)}}: $[K_{ij}, K_{kl}]=1$ for $(i,j)< (k,l)$, \\
\hspace{2em} {\bf{(P3)}}: $[K_{ik}, K_{ij} K_{kj}]=1$ \\
where in each relation, the indeces $i$, $j$, $k$, $l$ are distinct.
\end{thm}

\vspace{0.5em}

We construct a group extension of $\mathrm{P}\Sigma_n$ by $\mathfrak{S}_n$ as follows.
Let $B_n$ be Artin's braid group of $n$-strands. Artin \cite{Art} obtained the first finite presentation
\[\begin{split}
   B_n & = \langle \,\, \sigma_1, \ldots, \sigma_{n-1} \,|\, 
                     \sigma_{i+1} \sigma_i \sigma_{i+1} \sigma_i^{-1} \sigma_{i+1}^{-1} \sigma_i^{-1} \,\,\, (1 \leq i \leq n-2), \\
                  & \hspace{7.8em} [\sigma_i, \sigma_j] \,\,\, (1 \leq i<j \leq n-1, \,\, |i-j| \geq 2) \,\,\rangle.
  \end{split}\]
Moreover, Artin constructed a faithful representation of $B_n$ into $\mathrm{Aut}\,F_n$. For any $1 \leq i \leq n-1$,
the image of $\sigma_i$ by this representation is given by
\[ x_t \mapsto \begin{cases}
               x_i^{-1} x_{i+1} x_i & t=i, \\
               x_i & t=i+1, \\
               x_t & t \neq i, i+1.
              \end{cases}\]
We also denote this automorphism of $F_n$ by $\sigma_i$, and consider $B_n$ as a subgroup of $\mathrm{Aut}\,F_n$.
The subgroup $B_n \cap \mathrm{IA}_n$ of $B_n$ is the pure braid group of $n$-strands, and is denoted by $P_n$. (For details, see \cite{Bir} for example.)

\vspace{0.5em}

According to \cite{FRR}, we call
an element $\sigma \in \mathrm{Aut}\,F_n$ a braid-permutation automorphism of $F_n$ if there exists some $\mu \in \mathfrak{S}_n$ and some
$a_i \in F_n$ for any $1 \leq i \leq n$ such that
\[ x_i^{\sigma}=a_i^{-1} x_{\mu(i)} a_i \hspace{1em} (1 \leq i \leq n). \]
Let $\mathrm{BP}_n$ be the subgroup of $\mathrm{Aut}\,F_n$ consisting of all braid-permutation automorphisms of $F_n$.
The group $\mathrm{BP}_n$ is called the braid-permutation automorphism group of $F_n$.
The subgroup $\mathrm{BP}_n \cap \mathrm{IA}_n$ of $\mathrm{BP}_n$ is the basis-conjugating automorphism group $\mathrm{P}\Sigma_n$.
For any $1 \leq i \leq n-1$, let $s_i$ be the automorphism of $F_n$ which permutates $x_i$ and $x_{i+1}$, and fixes the other generators of the basis of $F_n$.
Clearly, the subgroup of $\mathrm{Aut}\,F_n$ generated by $s_i$ for all $1 \leq i \leq n-1$ is isomorphic to $\mathfrak{S}_n$.
Fenn-Rim\'{a}nyi-Rourke gave a finite presentation for $\mathrm{BP}_n$ as follows.
\begin{thm}[Fenn-Rim\'{a}nyi-Rourke \cite{FRR}]\label{T-FRR}
The group $BP_n$ generated by $\sigma_i$ and $s_i$ for $1 \leq i \leq n-1$ subject to relations: \\
\hspace{2em} {\bf{(B1)}}: $\sigma_i \sigma_{i+1} \sigma_i = \sigma_{i+1} \sigma_{i} \sigma_{i+1}$ for
                          $1 \leq i \leq n-2$, \\
\hspace{2em} {\bf{(B2)}}: $[\sigma_i, \sigma_j]=1$ for $|i-j| \geq 2$, \\
\hspace{2em} {\bf{(SY1)}}: $s_i^2=1$ for $1 \leq i \leq n-1$, \\
\hspace{2em} {\bf{(SY2)}}: $s_i s_{i+1} s_i = s_{i+1} s_i s_{i+1}$ for $1 \leq i \leq n-2$, \\
\hspace{2em} {\bf{(SY3)}}: $[s_i, s_j]=1$ for $|i-j| \geq 2$, \\
\hspace{2em} {\bf{(BP1)}}: $\sigma_i s_j = s_j \sigma_i$ for $|i-j| \geq 2$, \\
\hspace{2em} {\bf{(BP2)}}: $s_i s_{i+1} \sigma_i = \sigma_{i+1} s_i s_{i+1}$ for $1 \leq i \leq n-2$, \\
\hspace{2em} {\bf{(BP3)}}: $\sigma_i \sigma_{i+1} s_i = s_{i+1} \sigma_i \sigma_{i+1}$ for $1 \leq i \leq n-2$.
\end{thm}
\noindent
Then we have the following commutative diagram whose three rows are group extensions, and whose all vertical maps are injective.
\[\begin{CD}
1 @>>> \mathrm{IA}_n @>>> \mathrm{Aut}\,F_n @>{\rho}>> \mathrm{GL}(n,\Z) @>>> 1 \\
   @.     @AAA     @AAA      @AAA        @.  \\
1 @>>> \mathrm{P}\Sigma_n @>>> \mathrm{BP}_n @>{\rho|_{\mathrm{BP}_n}}>> \mathfrak{S}_n @>>> 1 \\
   @.     @AAA     @AAA      @AAA        @.  \\
1 @>>> P_n @>>> B_n @>{\rho|_{B_n}}>> \mathfrak{S}_n @>>> 1 \\
 \end{CD}\]
Here $\rho : \mathrm{Aut}\,F_n \rightarrow \mathrm{GL}(n,\Z)$ is the homomorphism induced from the abelianization of $F_n$.

\vspace{0.5em}

For any $k \geq 1$, set $\mathcal{P}_n(k):= \mathrm{P}\Sigma_n \cap \mathcal{A}_n(k)$. The filtration
$\mathrm{P}\Sigma_n = \mathcal{P}_n(1) \supset \mathcal{P}_n(2) \supset \cdots$ is a descending central filtration of $\mathrm{P}\Sigma_n$.
Set $\mathrm{gr}^k(\mathcal{P}_n):= \mathcal{P}_n(k)/\mathcal{P}_n(k+1)$.
In Proposition 3.5 in  \cite{S18}, Satoh showed that
for any $n \geq 2$ and $k \geq 1$, the image of
the restriction of
$\tilde{\tau}_k : \mathcal{A}_n(k) \rightarrow H^* \otimes_{\Z} \mathcal{L}_n(k+1)$ to $\mathcal{P}_n(k)$
is contained in $\mathfrak{p}_n(k)$.
\begin{defn}
We call the induced injective homomorphism
\[ \tau_k^P : \mathrm{gr}^k(\mathcal{P}_n) \rightarrow \mathfrak{p}_n(k) \]
from $\tilde{\tau}_k|_{\mathcal{P}_n(k)}$ the $k$-th Johnson homomorphism of $\mathrm{P}\Sigma_n$.
\end{defn}
\noindent
It is known that
for each $k \geq 1$, the Johnson homomorphism $\tau_k^P$ is $\mathrm{BP}_n/\mathrm{P}\Sigma_n = \mathfrak{S}_n$-equvariant.
and that $\tau_1^P$ induces the abelianization of $\mathrm{P}\Sigma_n$.
The graded sum $\mathrm{gr}(\mathcal{P}_n) := \oplus_{k \geq 1} \mathrm{gr}^k(\mathcal{P}_n)$ is considered as the Lie subalgebra of $\mathrm{gr}(\mathcal{A}_n)$,
and the graded sum $\tau^P := \oplus_{k \geq 1} \tau_k^P : \mathrm{gr}(\mathcal{P}_n) \rightarrow \mathfrak{p}_n$ is a Lie algebra homomorphism.

\vspace{0.5em}

Let $\mathrm{P}\Sigma_n=\mathrm{P}\Sigma_n{(1)} \supset \mathrm{P}\Sigma_n{(2)} \supset \cdots$ be the lower central series of $\mathrm{P}\Sigma_n$.
We have $\mathrm{P}\Sigma_n{(k)} \subset \mathcal{P}_n(k)$ for any $k \geq 1$.
Satoh \cite{S15} showed that $\mathrm{P}\Sigma_n{(k)} = \mathcal{P}_n(k)$ for $k \leq 3$.
In general, it is not known whether $\mathrm{P}\Sigma_n{(k)}$ coincides with $\mathcal{P}_n(k)$ or not for any $n \geq 3$ and $k \geq 4$.
Let $\mathrm{P}\Sigma_n^+$ be the subgroup of $\mathrm{P}\Sigma_n$ generaed by $K_{ij}$ for all $1 \leq j < i \leq n$.
We remark that
Darn\'{e} \cite{Da2} proved that $\mathrm{P}\Sigma_n^+ \cap \mathcal{A}_n(k)$ coincides with the $k$-the subgroup of the lower central series of $\mathrm{P}\Sigma_n^+$
for any $k \geq 1$.

\section{Trace maps}\label{S-trace}

Here we consider the trace maps and its variants to detect the cokernel of $\tau_k^P$ for any $k \geq 1$.
Let $\mathrm{Tr}_k^P : \mathfrak{p}_n(k) \rightarrow \mathcal{C}_n(k)$ be the restriction of $\mathrm{Tr}_k$ to $\mathfrak{p}_n(k)$.
In \cite{S03}, we showed that $\mathrm{Im}(\tau_k) \subset \mathrm{Ker}(\mathrm{Tr}_k)$ for any $n \geq 2$ and $k \geq 2$.
From this we see $\mathrm{Im}(\tau_k^P) \subset \mathrm{Ker}(\mathrm{Tr}_k^P)$.
Hence in order to study $\mathrm{Coker}(\tau_k^P)$,
it is interesting problem to determine the image of the trace map $\mathrm{Tr}_k^P$.

\vspace{0.5em}

It is easy to see that for any $k \geq 2$, $\langle x_1^k, \ldots, x_n^k \rangle \cap \mathrm{Im}(\mathrm{Tr}_k^P) = \{ 0 \}$.
(See also Proposition 3.7 in \cite{S18}.)
Set $\overline{\mathcal{C}}_n(k) := \mathcal{C}_n(k)/\langle x_1^k, \ldots, x_n^k \rangle$.
By compositing $\mathrm{Tr}_k^P$ and
the natural projection $\mathcal{C}_n(k) \rightarrow \overline{\mathcal{C}}_n(k)$,
we obtain the map $\overline{\mathrm{Tr}}_k^P : \mathfrak{p}_n(k) \rightarrow \overline{\mathcal{C}}_n(k)$.
As we will mention in Section {\rmfamily \ref{S-comp}}, the map $\overline{\mathrm{Tr}}_k^P$ is surjective for $k \leq 3$. Hence
we can describe the rational Johnson cokernel by using $\overline{\mathrm{Tr}}_k^P$ in these cases.
In general, however, $\overline{\mathrm{Tr}}_k^P$ is not surjective and it seems a difficult problem to determine its image.
Here we consider some reduced versions of $\overline{\mathrm{Tr}}_k^P$.
In \cite{S18}, one of such reduced versions of $\mathrm{Tr}_k^P$ was studied.
More precisely, consider the quotient $\mathfrak{S}_n$-module $\overline{S}^k(H):=S^k(H)/\langle x_1^k, \ldots, x_n^k \rangle$, and
the natural projection $\lambda_k : \mathcal{C}_n(k) \rightarrow \overline{S}^k(H)$.
In Proposition 3.7 in \cite{S18}, Satoh showed that the map $\lambda_k \circ \mathrm{Tr}_k^P$ is surjective for $n \geq 2$ and $k \geq 2$.

\vspace{0.5em}

Here we introduce the other trace maps which detect different comonents other than $\overline{S}^k(H_{\Q})$
in $\mathrm{Coker}(\tau_{k,\Q}^P)$.
\begin{defn}
For any $k \geq 2$,
let $N(k)$ be the $\mathfrak{S}_n$-invariant submodule of $\mathcal{C}_n(k)$ generated by
$x_{j_1} x_{j_2} \cdots x_{j_k}$ for all $1 \leq j_1, \ldots, j_k \leq n$ such that each $j_l$ $(1 \leq l \leq k)$ appears at least twice in
$j_1, j_2, \ldots, j_k$.
\end{defn}
\noindent
For example, basis of $N(k)$ for $2 \leq k \leq 4$ are given as follows.
\vspace{0.5em}
\begin{center}
{\renewcommand{\arraystretch}{1.3}
\begin{tabular}{|c|l|l|} \hline
  $k$  & basis of $N(k)$  &     \\ \hline
  $2$  & $x_i^2$        & $1 \leq i \leq n$  \\ \hline
  $3$  & $x_i^3$  & $1 \leq i \leq n$    \\ \hline 
  $4$  & $x_i^4$  & $1 \leq i \leq n$  \\
       & $x_i^2 x_j^2$, $x_ix_jx_ix_j$  & $1 \leq i<j \leq n$ \\ \hline
\end{tabular}}
\end{center}
\vspace{0.5em}
\noindent
Set $\widetilde{\mathcal{C}}_n(k) := \mathcal{C}_n(k)/N(k)$.
Let $\mu_k : \mathcal{C}_n(k) \rightarrow \widetilde{\mathcal{C}}_n(k)$ be the natural projection,
and set $\widetilde{\mathrm{Tr}}_k^P:= \mu_k \circ \mathrm{Tr}_k^P : \mathfrak{p}_n(k) \rightarrow \widetilde{\mathcal{C}}_n(k)$.
Then we have
\[ \mathrm{Im}(\tau_k^P) \subset \mathrm{Ker}(\mathrm{Tr}_k^P) = \mathrm{Ker}(\overline{\mathrm{Tr}}_k^P) \subset
   \mathrm{Ker}(\widetilde{\mathrm{Tr}}_k^P) \subset \mathfrak{p}_n(k). \]
\begin{thm}\label{T-main_1}
For any $n \geq 2$ and $k \geq 2$, the map $\widetilde{\mathrm{Tr}}_k^P$ is surjective.
\end{thm}
\textit{Proof.}
As a $\Z$-module, $\widetilde{\mathcal{C}}_n(k)$ is generated by elements
$x_{j_1} x_{j_2} \cdots x_{j_k} \in \widetilde{\mathcal{C}}_n(k)$
such that there exists some $1 \leq l \leq k$ such that $j_l \neq j_m$ for all $m \neq l$.
By using Lemma {\rmfamily \ref{L-S03_1}}, we have
\[\begin{split}
   \widetilde{\mathrm{Tr}}_k^P & (x_{j_l}^* \otimes [x_{j_l}, x_{j_{l+1}}, \ldots, x_{j_k}, x_{j_1}, \ldots, x_{j_{l-1}}, x_{j_l}]) \\
    & = x_{j_{l+1}} \cdots x_{j_k} x_{j_1} \cdots x_{j_l}
        - \mu_k \circ \varpi_k([x_{j_l}, x_{j_{l+1}}, \ldots, x_{j_k}, x_{j_1}, \ldots, x_{j_{l-1}}]).
  \end{split}\]
Since $\mathcal{L}_n(k) \subset \mathrm{Ker}(\varpi_k)$, we see that
\[ x_{j_1} x_{j_2} \cdots x_{j_k} = x_{j_{l+1}} \cdots x_{j_k} x_{j_1} \cdots x_{j_l} \in \mathrm{Im}(\widetilde{\mathrm{Tr}}_k^P). \]
Thus, $\widetilde{\mathrm{Tr}}_k^P$ is surjective.
$\square$

\begin{cor}\label{C-0}
For any $n \geq 2$ and $k \geq 2$, the module $\widetilde{\mathcal{C}}_n^{\Q}(k)$ injects into $\mathrm{Coker}(\tau_{k,\Q}^P)$
\end{cor}
\noindent
There is the natural surjective homomorphism $\widetilde{\mathcal{C}}_n(k) \rightarrow \overline{S}^k(H)$.
Hence in the case where $k=2$ and $3$, Proposition 3.7 in \cite{S18} is a corollary of Theorem {\rmfamily \ref{T-main_1}}.
\begin{cor}\label{C-1}
For any odd $k \geq 3$ and $n \geq k$, we have the natural surjective homomorphism $\widetilde{\mathcal{C}}_n(k) \rightarrow \wedge^k H$.
Therefore after tensoring with $\Q$, we see that the module $\wedge^k H_{\Q}$ injects into
$\mathrm{Coker}(\tau_{k,\Q}^P)$.
\end{cor}


\section{The Johnson cokernel in the case where $k \leq 4$}\label{S-comp}

In this section, we determine the image and the cokernel of $\tau_k^P$ for $1 \leq k \leq 4$ by using the trace maps considered
in Section {\rmfamily \ref{S-trace}}. First, we recall some results in \cite{S15}.

\vspace{0.5em}

For $k=1$, the first Johnson homomorphism $\tau_1^P : \mathrm{gr}^1(\mathcal{P}_n) \rightarrow \mathfrak{p}_n(1)$ is an isomorphism, and
$\mathrm{Im}(\tau_1^P)$ is a free abelian group of rank $n(n-1)$ with a basis
$\{x_i^* \otimes [x_j, x_i] \,|\, 1 \leq i \neq j \leq n \}$.
From Lemma 2.3 in \cite{S15}, we see that $\mathrm{Im}(\tau_2^P)$ is a free abelian group of rank $\frac{n(n-1)^2}{2}$ with a basis
\begin{equation}\begin{split}
    \{ x_i^* & \otimes [x_{j_1}, x_{j_2}, x_i] \,|\, 1 \leq j_2 < j_1 \leq n, \hspace{0.5em} i \neq j_1, j_2 \} \\
   & \cup \{ x_i^* \otimes [x_i, x_j, x_i] - x_j^*  \otimes [x_j, x_i, x_j] \,|\, 1 \leq i \neq j \leq n \}
  \end{split} \label{eq-basis-tau_2}
\end{equation}
for $n \geq 3$.
\begin{pro}\label{P-1}
For $n \geq 3$, we have \\
$(1)$\,\,$\overline{\mathrm{Tr}}_2^P=\widetilde{\mathrm{Tr}}_2^P$ is surjective. \\
$(2)$\,\,$\mathrm{Im}(\tau_2^P)=\mathrm{Ker}(\overline{\mathrm{Tr}}_2^P)$.
\end{pro}
\textit{Proof.}
(1)\,\,By definition, we see $\overline{\mathrm{Tr}}_2^P=\widetilde{\mathrm{Tr}}_2^P$.
For any generator $x_i x_j$ ($i \neq j$) of $\overline{\mathcal{C}}_n(2)$, we have
\[ \overline{\mathrm{Tr}}_2^P(x_i^* \otimes [x_i, x_j, x_i]) = x_jx_i=x_ix_j. \]
Hence $\overline{\mathrm{Tr}}_2^P$ is surjective. \\
(2)\,\,It suffices to show $\mathrm{Im}(\tau_2^P) \supset \mathrm{Ker}(\overline{\mathrm{Tr}}_2^P)$.
Take any
\[ x := \sum_{\substack{1 \leq i \leq n\\[1pt] 1 \leq j_2 < j_1 \leq n}} a_{j_1, j_2}^i x_i^* \otimes [x_{j_1}, x_{j_2}, x_i] \in \mathrm{Ker}(\overline{\mathrm{Tr}}_2^P). \]
Since $\mathrm{Im}(\tau_2^P) \subset \mathrm{Ker}(\overline{\mathrm{Tr}}_2^P)$, and from (\ref{eq-basis-tau_2}),
we may assume
\[ x=\sum_{1 \leq i < j \leq n} a_{i, j}^i x_i^* \otimes [x_i, x_j, x_i]. \]
Then we obtain
\[ 0=\tau_2^P(x)=\sum_{1 \leq i < j \leq n} a_{i, j}^i x_ix_j, \]
and hence $a_{i, j}^i=0$ for any $i<j$. Thus $x=0$. $\square$

\begin{cor}\label{C-3}
For $n \geq 3$,
\[ 0 \rightarrow \mathrm{gr}^2(\mathcal{P}_n) \xrightarrow{\tau_2^P} \mathfrak{p}_n(2) \xrightarrow{\overline{\mathrm{Tr}}_2^P} \overline{\mathcal{C}}_n(2) \rightarrow 0 \]
is $\mathfrak{S}_n$-equivariant exact.
\end{cor}

Next, we consider the case where $k=3$. We prepare the following lemma.
\begin{lem}\label{L-1-1}
For any $k \geq 3$ and any $1 \leq i \leq n$, if $1 \leq j_1, \ldots, j_k \leq n$ and $i \neq j_1, \ldots, j_k$, then
$x_i^* \otimes [x_{j_1}, \ldots, x_{j_k}, x_i] \in \mathrm{Im}(\tau_k^P)$.
\end{lem}
\textit{Proof.}
We prove this lemma by induction on $k \geq 2$. For $k=2$, it is clear from (\ref{eq-basis-tau_2}).
Assume $k \geq 3$. By the inductive hypothesis, we see
\[ x_i^* \otimes [x_{j_1}, \ldots, x_{j_k-1}, x_i] \in \mathrm{Im}(\tau_{k-1}^P). \]
If we write $x_i^* \otimes [x_{j_1}, \ldots, x_{j_k-1}, x_i]$ as $\tau_{k-1}^P(\sigma)$ for $\sigma \in \mathcal{P}_n(k-1)$,
then we have
\[\begin{split}
   \tau_k^P([\sigma, K_{i, j_k}]) &= [\tau_{k-1}^P(\sigma), \tau_1^P(K_{i, j_k})]  \\
    &=[x_i^* \otimes [x_{j_1}, \ldots, x_{j_k-1}, x_i], x_i^* \otimes [x_i, x_{j_k}] ] \\
    &=x_i^* \otimes [x_{j_1}, \ldots, x_{j_k-1}, [x_i, x_{j_k}] ] - x_i^* \otimes [x_{j_1}, \ldots, x_{j_k-1}, x_i, x_{j_k}]  \\
    &= - x_i^* \otimes [x_{j_1}, \ldots, x_{j_k}, x_i].
 \end{split}\]
At the last equation, we use the Jacobi identity. Hence the induction proceeds.
$\square$

\begin{pro}\label{P-2}
For $n \geq 3$, we have\\
$(1)$\,\,$\overline{\mathrm{Tr}}_3^P=\widetilde{\mathrm{Tr}}_3^P$ is surjective. \\
$(2)$\,\,$\mathrm{Im}(\tau_3^P)=\mathrm{Ker}(\overline{\mathrm{Tr}}_3^P)$.\\
$(3)$\,\,$\mathrm{Im}(\tau_3^P)$ is a free abelian group with a basis
\begin{equation}\begin{split}
    \{ x_i^* & \otimes [x_{j_1}, x_{j_2}, x_{j_3}, x_i] \,|\, j_1 > j_2 \leq j_3, \hspace{0.5em} i \neq j_1, j_2, j_3 \} \\
   & \cup \{ x_i^* \otimes [x_i, x_j, x_i, x_i] - x_j^* \otimes [x_j, x_i, x_i, x_j] \,|\, 1 \leq i \neq j \leq n \} \\
   & \cup \{ x_i^* \otimes [x_i, x_j, x_l, x_i] - x_j^* \otimes [x_j, x_l, x_i, x_j] \,|\, l \neq i > j \neq l \}
  \end{split}
\end{equation}
\end{pro}
\textit{Proof.}
(1)\,\,By definition, we see $\overline{\mathrm{Tr}}_3^P=\widetilde{\mathrm{Tr}}_3^P$.
$\overline{\mathcal{C}}_n(3)$ is generated by
$x_ix_jx_l$ and $x_i^2x_j$
for all distinct $1 \leq i, j, l \leq n$. Since
\[\begin{split}
   \overline{\mathrm{Tr}}_3^P(x_i^* \otimes [x_i, x_j, x_l, x_i]) & = x_jx_lx_i=x_ix_jx_l, \\
   \overline{\mathrm{Tr}}_3^P(x_i^* \otimes [x_i, x_j, x_i, x_i]) & = x_jx_i^2=x_i^2x_j, \\
  \end{split}\]
we see that $\overline{\mathrm{Tr}}_3^P$ is surjective.\\
(2)\,\,It suffices to show $\mathrm{Im}(\tau_3^P) \supset \mathrm{Ker}(\overline{\mathrm{Tr}}_3^P)$.
Take any
\[ x := \sum_{\substack{1 \leq i \leq n\\[1pt] j_1>j_2 \leq j_3}} a_{j_1, j_2, j_3}^i x_i^* \otimes [x_{j_1}, x_{j_2}, x_{j_3}, x_i] \in \mathrm{Ker}(\overline{\mathrm{Tr}}_3^P). \]
Observe
\begin{equation}\begin{split}
   [x_{j_1}, x_{j_2}, x_i, x_i] &= -[x_i, [x_{j_1}, x_{j_2}], x_i] = -[x_i, x_{j_1}, x_{j_2}, x_i] + [x_i, x_{j_2}, x_{j_1}, x_i], \\
   \tau_3^P([K_{ij},K_{ji}, K_{ji}]) &= x_i^* \otimes [x_i, x_j, x_i, x_i] - x_j^* \otimes [x_j, x_i, x_i, x_j]
  \end{split} \label{eq-basis-tau_3}
\end{equation}
for any $1 \leq i, j, j_1, j_2 \leq n$. Remark that
\[ [x_j,x_i,x_j,x_i] - [x_j,x_i,x_i,x_j]=[x_j,x_i,[x_j,x_i]]=0 \]
by the Jacobi identity.
Since $\mathrm{Im}(\tau_3^P) \subset \mathrm{Ker}(\overline{\mathrm{Tr}}_3^P)$, and from Lemma {\rmfamily \ref{L-1-1}} and (\ref{eq-basis-tau_3}),
we may assume
\[ x=\sum_{\substack{1 \leq i, j, l \leq n\\[1pt] j, l \neq i}} a_{i, j, l}^i x_i^* \otimes [x_i, x_j, x_l, x_i]. \]
Furthermore, by observing
\[ \tau_3^P([[K_{ij},K_{il}], K_{ji}]) = x_j^* \otimes [x_j, x_l, x_i, x_j] - x_i^* \otimes [x_i, x_j, x_l, x_i] \]
for any distinct $i, j$ and $l$, we may also assume
\[ x=\sum_{\substack{i< j, l\\[1pt] j \neq l}} a_{i, j, l}^i x_i^* \otimes [x_i, x_j, x_l, x_i] + \sum_{1 \leq i \neq j \leq n} a_{i, j, j}^i x_i^* \otimes [x_i, x_j, x_j, x_i]. \]
Then we obtain
\[ 0=\tau_3^P(x)=\sum_{\substack{i< j, l\\[1pt] j \neq l}} a_{i, j, l}^i x_ix_jx_l + \sum_{1 \leq i \neq j \leq n} a_{i, j, j}^i x_ix_j^2, \]
and hence $a_{i, j, l}^i=0$ for any $i, j$ and $l$. Thus $x=0$. \\
(3)\,\,From Part (2), we see that $\mathrm{Im}(\tau_3^P)$ is generated by
\[\begin{cases}
   x_i^* \otimes [x_{j_1}, x_{j_2}, x_{j_3}, x_i] \hspace{1em} & i \neq j_1, j_2, j_3, \\
   x_i^* \otimes [x_i, x_j, x_i, x_i] - x_j^* \otimes [x_j, x_i, x_i, x_j] & i \neq j, \\
   x_i^* \otimes [x_i, x_j, x_l, x_i] - x_j^* \otimes [x_j, x_l, x_i, x_j] & i \neq j \neq l \neq i
\end{cases}\]
as a $\Z$-module.
$\square$

\begin{cor}\label{C-4}
For $n \geq 3$,
\[ 0 \rightarrow \mathrm{gr}^3(\mathcal{P}_n) \xrightarrow{\tau_3^P} \mathfrak{p}_n(3) \xrightarrow{\overline{\mathrm{Tr}}_3^P} \overline{\mathcal{C}}_n(3) \rightarrow 0 \]
is $\mathfrak{S}_n$-equivariant exact.
\end{cor}
\noindent
In \cite{ES1}, we showed that $\mathcal{C}_n^{\Q}(3)$ is isomorphic to $S^3 H_{\Q} \oplus \wedge^3 H_{\Q}$. From this fact and Corollary {\rmfamily \ref{C-4}}, we have the
$\mathfrak{S}_n$-equivariant exact sequence
\[ 0 \rightarrow \mathrm{gr}^3(\mathcal{P}_n) \xrightarrow{\tau_3^P} \mathfrak{p}_n(3) \xrightarrow{} \overline{S}^3 H_{\Q} \oplus \wedge^3 H_{\Q} \rightarrow 0 \]
for $n \geq 3$.

\vspace{0.5em}

Finally, we consider the case where $k=4$.
At first, we prepare a few lemmas.
\begin{lem}\label{L-1-2}
For any $n \geq 3$, the module $\mathfrak{p}_n(4)$ is generated by all
\begin{equation}\begin{cases}
   x_i^* \otimes [x_{j_1}, x_{j_2}, x_{j_3}, x_{j_4}, x_i] \hspace{1em} & i \neq j_1, j_2, j_3, j_4, \\
   x_i^* \otimes [x_i, x_{j_2}, x_{j_3}, x_{j_4}, x_i] & i \neq j_2, j_3, j_4, \\
   x_i^* \otimes [x_i, x_{j_2}, x_{i}, x_{j_4}, x_i] & i \neq j_2, j_4, \\
   x_i^* \otimes [x_i, x_{j_2}, x_{j_3}, x_i, x_i] & i \neq j_2, j_3, \\
   x_i^* \otimes [x_i, x_{j_2}, x_i, x_i, x_i] & i \neq j_2 \\
  \end{cases} \label{eq-L-1-2}
\end{equation}
for $1 \leq i, j_1, j_2, j_3, j_4 \leq n$ as a $\Z$-module.
\end{lem}
\textit{Proof.}
By the definition, $\mathfrak{p}_n(4)$ is generated by elements $x_i^* \otimes [x_{j_1}, x_{j_2}, x_{j_3}, x_{j_4}, x_i]$
for all $1 \leq i, j_1, j_2, j_3, j_4 \leq n$ as a $\Z$-module. Set $\omega := x_i^* \otimes [x_{j_1}, x_{j_2}, x_{j_3}, x_{j_4}, x_i]$.
If $j_1, \ldots, j_4 \neq i$, then $\omega$ is an element of first type.
If $j_1, j_2 \neq i$ and $j_3=i$, then by
\[\begin{split}
   [x_{j_1}, x_{j_2}, x_i, x_{j_4}, x_i] &= -[x_i, [x_{j_1}, x_{j_2}], x_{j_4}, x_i] \\
    & = -[x_i, x_{j_1}, x_{j_2}, x_{j_4}, x_i] + [x_i, x_{j_2}, x_{j_1}, x_{j_4}, x_i]
  \end{split}\]
from the Jacobi identity, we see that $\omega$ can be written as a sum of elements the form in (\ref{eq-L-1-2}).
If $j_1, j_2, j_3 \neq i$ and $j_4=i$, then we can show by an argument similar to that in case where $j_1, j_2 \neq i$ and $j_3=i$.
If $j_1=i$ and $j_2 \neq i$, or if $j_2=i$ and $j_1 \neq i$, then we see that $\omega$ is an element of the above form.
$\square$

\begin{lem}\label{L-1-2}
For any $n \geq 3$, the map $\overline{\mathrm{Tr}}_4^P$ is not surjective, and its cokernel is a free abelian group of rank $n(n-1)/2$.
\end{lem}
\textit{Proof.}
A basis of $\overline{\mathcal{C}}_n(4)$ is given by as follows:
\begin{itemize}
\item $x_ax_bx_cx_d$ for any distinct $a, b, c, d$ such that $a<b, c, d$,
\item $x_a^2 x_bx_c$ for any distinct $a, b, c$,
\item $x_a x_b x_a x_c$ for any distinct $a, b, c$ such that $b<c$,
\item $x_a^2 x_b^2$ for any $a < b$,
\item $x_ax_bx_ax_b$ for any $a < b$,
\item $x_a^3 x_b$ for any $a \neq b$.
\end{itemize}
On the other hand, the images of the generators of $\mathfrak{p}_n(4)$ shown in (\ref{eq-L-1-2}) by $\overline{\mathrm{Tr}}_4^P$ are given by
\[\begin{cases}
   \overline{\mathrm{Tr}}_4^P(x_i^* \otimes [x_{j_1}, x_{j_2}, x_{j_3}, x_{j_4}, x_i])=0 \hspace{1em} & i \neq j_1, j_2, j_3, j_4, \\
   \overline{\mathrm{Tr}}_4^P(x_i^* \otimes [x_i, x_{j_2}, x_{j_3}, x_{j_4}, x_i])= x_{j_2}x_{j_3}x_{j_4}x_i & i \neq j_2, j_3, j_4, \\
   \overline{\mathrm{Tr}}_4^P(x_i^* \otimes [x_i, x_{j_2}, x_{i}, x_{j_4}, x_i])=2x_{j_2}x_ix_{j_4}x_i - x_{j_2}x_{j_4}x_i^2 & i \neq j_2, j_4, \\
   \overline{\mathrm{Tr}}_4^P(x_i^* \otimes [x_i, x_{j_2}, x_{j_3}, x_i, x_i])=2x_{j_2}x_ix_{j_3}x_i - x_{j_2}x_{j_3}x_i^2 & i \neq j_2, j_3, \\
   \overline{\mathrm{Tr}}_4^P(x_i^* \otimes [x_i, x_{j_2}, x_i, x_i, x_i])=x_{j_2}x_i^3 & i \neq j_2. \\
  \end{cases}\]
Hence the elements of the basis types of $x_ax_bx_cx_d$, $x_a^2 x_bx_c$, $x_a x_b x_a x_c$ and $x_a^3 x_b$ are in $\mathrm{Im}(\overline{\mathrm{Tr}}_4^P)$.
Indeed, for example,
\[ \overline{\mathrm{Tr}}_4^P(x_c^* \otimes [x_c, x_a, x_a, x_b, x_c])= x_a^2 x_bx_c. \]
Therefore we see that $\mathrm{Coker}(\overline{\mathrm{Tr}}_4^P)$ is isomorphic to the quotient module of the free $\Z$-module with basis $x_a^2 x_b^2$
and $x_ax_bx_ax_b$ for all $a < b$ divided by the submodule generated by $2x_ax_bx_ax_b - x_a^2x_b^2$ for all $a < b$. This completes the proof. $\square$

\begin{pro}\label{P-3}
For $n \geq 3$, we have $\mathrm{Im}(\tau_4^P)=\mathrm{Ker}(\overline{\mathrm{Tr}}_4^P)$.
\end{pro}
\textit{Proof.}
It suffices to show $\mathrm{Im}(\tau_4^P) \supset \mathrm{Ker}(\overline{\mathrm{Tr}}_4^P)$.
Take any $x \in \mathrm{Ker}(\overline{\mathrm{Tr}}_4^P)$. 
For distinct $1 \leq i, j, l \leq n$, set
\[ D_1 := x_i^* \otimes [x_j, x_l, x_i, x_i] + x_j^* \otimes [x_j, x_l, x_i, x_j] - x_j^* \otimes [x_j, x_i, x_l, x_j] \in \mathfrak{p}_n(3). \]
Since $D_1 \in \mathrm{Ker}(\overline{\mathrm{Tr}}_3^P) = \mathrm{Im}(\tau_3^P)$, there exists some $L_1 \in \mathcal{P}_n(3)$ such that $\tau_3^P(L_1)=D_1$.
From Lemma {\rmfamily \ref{L-1-1}} and
\begin{equation}\begin{split}
   \tau_4^P([K_{ij},K_{ji}, K_{ji}, K_{ji}]) &= x_i^* \otimes [x_i, x_j, x_i, x_i, x_i] - x_j^* \otimes [x_j, x_i, x_i, x_i, x_j], \\
   \tau_4^P([K_{ij},K_{ji}, K_{ji}, K_{ij}]) &= x_i^* \otimes [x_i, x_j, x_j, x_i, x_i] - x_j^* \otimes [x_j, x_i, x_j, x_i, x_j], \\
   \tau_4^P([K_{ji}, L_1]) &= x_i^* \otimes [x_i, x_j, x_l, x_i, x_i] - 2 x_j^* \otimes [x_j, x_i, x_l, x_i, x_j] \\
                              & \hspace{2em} + x_j^* \otimes [x_j, x_i, x_i, x_l, x_j], \\
  \end{split} \label{eq-P-3-1}
\end{equation}
for any distinct $1 \leq i, j, l \leq n$, by reducing the generators shown in (\ref{eq-L-1-2}) by using these elements
we may assume
\[\begin{split}
  x=\sum_{\substack{1 \leq i, j, l, m \leq n\\[1pt] j, l, m \neq i}} \alpha_{j, l, m}^i x_i^* \otimes [x_i, x_j, x_l, x_m, x_i]
            + \sum_{\substack{1 \leq i, j, m \leq n\\[1pt] j, m \neq i}} \beta_{j, m}^i x_i^* \otimes [x_i, x_j, x_i, x_m, x_i].
  \end{split}\]

For distinct $1 \leq i, j, l \leq n$, set
\[ D_2 := x_i^* \otimes [x_i, x_j, x_l, x_i] - x_l^* \otimes [x_l, x_i, x_j, x_l] \in \mathfrak{p}_n(3). \]
Since $D_2 \in \mathrm{Ker}(\overline{\mathrm{Tr}}_3^P) = \mathrm{Im}(\tau_3^P)$, there exists some $L_2 \in \mathcal{P}_n(3)$ such that $\tau_3^P(L_2)=D_2$.
Thus we have
\[\begin{split}
    \tau_4^P([L_2, K_{ji}]) &= x_i^* \otimes [x_i, x_j, x_i, x_l, x_i] - x_l^* \otimes [x_l, x_i, x_j, x_i, x_l] \\
                              & \hspace{2em} + x_l^* \otimes [x_l, x_i, x_i, x_j, x_l] - x_j^* \otimes [x_j, x_i, x_l, x_i, x_j].
  \end{split}\]
Moreover, by the Jacobi identity, we have
\[ [x_i, x_j, x_i, x_j, x_i] =- [x_i, x_j, [x_i,x_j], x_i] -[x_j,[x_i,x_j], x_i,x_i]=[x_i,x_j,x_j,x_i,x_i], \]
and hence from the second equation of (\ref{eq-P-3-1}), we may assume $\beta_{j, m}^i=0$ if $j \neq m$, or $j=m$ and $i>j$.

\vspace{0.5em}

On the other hand, consider the following equations.
\begin{equation}
 \tau_4^P([[K_{ac},K_{ad}], [K_{ca},K_{cb}]]) = x_a^* \otimes [x_a, x_b, x_c, x_d, x_a] - x_c^* \otimes [x_c, x_d, x_a, x_b, x_c] \label{eq-P-3-2}
\end{equation}
for $1 \leq a, b, c, d \leq n$ such that $c \neq d$ and $b \neq c$.
\begin{equation}
 \tau_4^P([K_{ab}, K_{ac}, K_{ad}, K_{ba}]) = x_a^* \otimes [x_a, x_b, x_c, x_d, x_a] - x_b^* \otimes [x_b, x_c, x_d, x_a, x_b] \label{eq-P-3-3}
\end{equation}
for $1 \leq a, b, c, d \leq n$ such that $b \neq c$.
\begin{equation}
 \tau_4^P([K_{ad}, [K_{da}, K_{db}, K_{dc}]]) = x_i^* \otimes [x_a, x_b, x_c, x_d, x_a] - x_d^* \otimes [x_d, x_a, x_b, x_c, x_d] \label{eq-P-3-4}
\end{equation}
for $1 \leq a, b, c, d \leq n$ such that $b, c \neq d$.
Then from (\ref{eq-P-3-4}) for the case where $b=c$, we see that we may assume $\alpha_{j,l,l}^i=0$.
Similarly, from (\ref{eq-P-3-2}) for the case where $b=d$, we see that we may assume $\alpha_{j,l,j}^i=0$ if $i>j$.
Furtheremore, for distinct $i, j, l, m$, by using (\ref{eq-P-3-2}), (\ref{eq-P-3-3}) and (\ref{eq-P-3-4}), we may assume $\alpha_{j, l, m}^i=0$
for the case other than $i< j, l, m$.
Therefore, $x$ can be written as
\[\begin{split}
  x &=\sum_{\substack{i< j, l, m \\[1pt] j, l, m : \text{distinct}}} \alpha_{j, l, m}^i x_i^* \otimes [x_i, x_j, x_l, x_m, x_i]
      + \sum_{\substack{i< l \\[1pt] j \neq i, l}} \alpha_{j, l, j}^i x_i^* \otimes [x_i, x_j, x_l, x_j, x_i] \\
    & + \sum_{i, j, l : \text{distinct}} \alpha_{j, j, l}^i x_i^* \otimes [x_i, x_j, x_j, x_l, x_i]
      + \sum_{i \neq j} \alpha_{j, j, j}^i x_i^* \otimes [x_i, x_j, x_j, x_j, x_i] \\
    &+ \sum_{i<j} \beta_{j, j}^i x_i^* \otimes [x_i, x_j, x_i, x_j, x_i].
  \end{split}\]
By observing
\[\begin{split}
  0=\overline{\mathrm{Tr}}_4^P(x) &=\sum_{\substack{i< j, l, m \\[1pt] j, l, m : \text{distinct}}} \alpha_{j, l, m}^i x_j x_l x_m x_i
      + \sum_{\substack{i< l \\[1pt] j \neq i, l}} \alpha_{j, l, j}^i x_j x_l x_j x_i \\
    & + \sum_{i, j, l : \text{distinct}} \alpha_{j, j, l}^i x_j^2 x_l x_i
      + \sum_{i \neq j} \alpha_{j, j, j}^i x_j^3 x_i \\
    &+ \sum_{i<j} \beta_{j, j}^i (2x_i x_j x_i x_j - x_i^2x_j^2),
  \end{split}\]
we obtain that all of the coefficients in the above equation is $0$, and hence $x=0$.
$\square$

\begin{cor}\label{C-5}
For $n \geq 3$,
\[ 0 \rightarrow \mathrm{gr}^4(\mathcal{P}_n) \xrightarrow{\tau_4^P} \mathfrak{p}_n(4) \xrightarrow{\overline{\mathrm{Tr}}_4^P} \overline{\mathcal{C}}_n(4)
    \rightarrow Q \rightarrow 0 \]
is $\mathfrak{S}_n$-equivariant exact. In particular, $\mathrm{Im}(\tau_4^P)$ is a free abelian group of rank $\frac{1}{4}n^2(n-1)^2(n+1)$.
\end{cor}
\noindent
We remark that for any $1 \leq a < b \leq n$, since
\[ 0 \neq 2x_ax_bx_ax_b - x_a^2x_b^2 \in \mathrm{Im}(\mathrm{Tr}_4^P) \cap N(4) \]
from the above calculation, we see $\mathrm{Ker}(\overline{\mathrm{Tr}}_4^P) \subsetneq \mathrm{Ker}(\widetilde{\mathrm{Tr}}_4^P)$. In general, we have the following.
\begin{lem}\label{L-add-1}
For any $k \geq 4$, $\mathrm{Ker}(\overline{\mathrm{Tr}}_k^P) \subsetneq \mathrm{Ker}(\widetilde{\mathrm{Tr}}_k^P)$.
Namely, we cannot detect all $\mathfrak{S}_n$-components in $\mathrm{Coker}(\tau_{k,\Q}^P)$ by $\widetilde{\mathrm{Tr}}_k^P$ for any $k \geq 4$.
\end{lem}
\textit{Proof.}
For any distinct $1 \leq i, j \leq n$, cosider an element $x_i^* \otimes [x_i, x_j, x_i, x_j, \ldots, x_j, x_i] \in \mathfrak{p}_n(k)$.
Then we have
\[ \mathrm{Tr}_k^P(x_i^* \otimes [x_i, x_j, x_i, x_j, \ldots, x_j, x_i]) = 2x_jx_ix_j^{k-3}x_i - x_i^2x_j^{k-2} \in N(k). \]
This shows $\mathrm{Ker}(\overline{\mathrm{Tr}}_k^P) \subsetneq \mathrm{Ker}(\widetilde{\mathrm{Tr}}_k^P)$.
$\square$

\vspace{0.5em}

On the other hand, for $k=4$, the $\mathfrak{S}_n$-irreducible component in $\mathrm{Coker}(\tau_4^P)$
which cannot detect $\widetilde{\mathrm{Tr}}_4^P$ is
the module
\[ A := \langle x_a^* \otimes [x_a,x_b,x_a,x_b,x_a] \,|\, 1 \leq a<b \leq n \rangle. \]
This component can be detected as follows. Let $J$ be the quotient module of $\wedge^2 H \otimes \wedge^2 H$ by the submodule generated by
\begin{itemize}
\item $(v \wedge w) \otimes (x \wedge y) - (x \wedge y) \otimes (v \wedge w)$,
\item $(v \wedge w) \otimes (x \wedge y) - (x \wedge w) \otimes (v \wedge y) - (v \wedge x) \otimes (w \wedge y)$
\end{itemize}
for any $v, w, x, y \in H$. Then $J$ is a free abelian group of rank $n^2(n^2-1)/12$.
As a $\mathrm{GL}(n,\Q)$-module, $J_{\Q}$ is irreducible and isomorphic to $[2,2]$.
In \cite{S09}, we construct a trace map which detect $J_{\Q}$ in the cokernel of $\tau_{4,\Q}$.
For $i=1, 2$, let $f_i : H^{\otimes 4} \rightarrow J$ be the projection defined by
\[ f_i(v \otimes w \otimes x \otimes y) = \begin{cases}
                                            (v \wedge x) \cdot (w \wedge y), \hspace{1em} & i=1, \\
                                            (v \wedge y) \cdot (w \wedge x), \hspace{1em} & i=2. \\
                                          \end{cases}\]
Set
\[ \mathrm{Tr}_J := f_1 \circ \Phi^4 - 2 (f_2 \circ \Phi^4) : H^* \otimes_{\Z} H^{\otimes 5} \rightarrow J. \]
By Proposition 4.2 in \cite{S09} and Darn\'{e}'s result in \cite{Dar}, we have 
$\mathrm{Im}(\tau_4) \subset \mathrm{Ker}(\mathrm{Tr}_J)$. Denote by $\mathrm{Tr}_J^P : \mathfrak{p}_n(4) \rightarrow J$
the restriction of $\mathrm{Tr}_J$ to $\mathfrak{p}_n(4)$.
Since we see
\[ \mathrm{Tr}_J(x_a^* \otimes [x_a,x_b,x_a,x_b,x_a]) = 5(x_a \wedge x_b) \cdot (x_a \wedge x_b) \neq 0, \]
we see that $\mathrm{Tr}_J$ detects $A$.

\vspace{0.5em}

Here we give the list of the rank of the free abelian groups $\mathrm{gr}^k(\mathcal{P}_n)$, $\mathfrak{p}_n(k)$, $\mathcal{C}_n(k)$ and $\mathrm{Coker}(\overline{\mathrm{Tr}}_k^P)$
for $n \geq 3$ and $1 \leq k \leq 4$.
\vspace{0.5em}
\begin{center}
{\renewcommand{\arraystretch}{1.1}
\begin{tabular}{|c|l|l|l|l|} \hline
  $k$  & $\mathrm{gr}^k(\mathcal{P}_n)$ & $\mathfrak{p}_n(k)$ & $\overline{\mathcal{C}}_n(k)$ & $\mathrm{Coker}(\overline{\mathrm{Tr}}_k^P)$ \\ \hline
  $1$  & $n(n-1)$                       & $n(n-1)$            & $0$                 & $0$  \\ \hline
  $2$  & $n(n-1)^2/2$                   & $n^2(n-1)/2$        & $n(n-1)/2$          & $0$  \\ \hline
  $3$  & $n(n-1)^2(n+1)/3$              & $n^2(n^2-1)/3$      & $n(n^2-1)/3$        & $0$  \\ \hline
  $4$  & $n^2(n-1)^2(n+1)/4$            & $n^3(n^2-1)/4$      & $n(n-1)(n^2+n-2)/4$ & $n(n-1)/2$  \\ \hline
\end{tabular}}
\end{center}

\section{Observations of the Johnson images and cokernels for $k \geq 5$}\label{S-OBS}

In this section, we give some observations of the dimensions of the Johnson images and cokernels for the case of $k \geq 5$ based on some direct computer calculations.
We will use the following notations. For a composition $\alpha=(\alpha_1,\alpha_2, \ldots ,\alpha_n)$ of $k$,
let $p(\alpha)$ be the set of all permutations of $1^{\alpha_1}2^{\alpha_2} \cdots n^{\alpha_n}$.
We denote $\mathcal{L}_n(k,\alpha)$ by the subspace of $\mathcal{L}_n(k)$ generated by all left-normed elements $[x_{i_1},x_{i_2}, \ldots ,x_{i_k}]$ such that $(i_1,i_2, \ldots ,i_k) \in p(\alpha)$.
Notice that $\alpha_j=\sharp\{1 \leq u \leq k \ | \ i_u=j\} \ (1 \leq j \leq n)$. Similarly, let $\overline{\mathcal{C}}_n(k,\alpha)$ be the subspace of $\overline{\mathcal{C}}_n(k)$ generated by all cosets (necklace) of $x_{i_1}x_{i_2} \cdots x_{i_k}$ for $(i_1,i_2, \ldots ,i_k) \in p(\alpha)$. \\
\quad For a $[x_{i_1}, \ldots ,x_{i_k}] \in \mathcal{L}_n(k,\alpha)$, 
through the natural embedding $\mathcal{L}_n(k)$ into $H^{\otimes k}$,
we consider $[x_{i_1}, \ldots ,x_{i_k}] \in H^{\otimes k}$, and
take the expression $[x_{i_1}, \ldots ,x_{i_k}]=\sum_{j=1}^{n}x_j \otimes b_j(i_1,\ldots ,i_k)$, where $b_j(i_1,\ldots ,i_k) \in H^{\otimes {k-1}}$. Since $\varpi^k([x_{i_1}, \ldots ,x_{i_k}])=0$, we have 
\begin{align*}
&\overline{\mathrm{Tr}}_k^P(x_i^* \otimes [x_{i_1}, \ldots ,x_{i_k},x_i])=
\overline{\mathrm{Tr}}_k^P(x_i^* \otimes 
([x_{i_1}, \ldots ,x_{i_k}] \otimes x_i-x_i \otimes [x_{i_1}, \ldots ,x_{i_k}]))\\
&
=\varpi^k(b_i(i_1, \ldots ,i_k) \otimes x_i-[x_{i_1},\ldots ,x_{i_k}])
=\varpi^k(b_i(i_1, \ldots ,i_k) \otimes x_i) \in  \overline{\mathcal{C}}_n(k,\alpha).
\end{align*}
Moreover 
\[
c_\alpha=\dim_{\Q} \left\langle
\varpi^k(b_i(i_1, \ldots ,i_k) \otimes x_i) \ | \ 1 \leq i \leq n,(i_1,\ldots ,i_k) \in p(\alpha)
\right\rangle_{\Q}
\]
is invariant under the left (variable permutation) action of $\mathfrak{S}_n$ on $\alpha$. Thus we may consider the case where
$\alpha=(\alpha_1 \geq \alpha_2 \geq \cdots \geq \alpha_n)$ is a partion of $k$. We obtain 
\[
\dim_{\Q}\mathrm{Im}(\overline{\mathrm{Tr}}_{k,\Q}^P)=
\sum_{\alpha\vdash k}|\mathfrak{S}_n \cdot \alpha| \cdot c_\alpha,
\]
where $|\mathfrak{S}_n \cdot \alpha|$ is the order of $\mathfrak{S}_n$-orbit of $\alpha$. Since the surjectivity of $\widetilde{\mathrm{Tr}}_k^P$ (Theorem \ref{T-main_1}), if $\alpha$ contains $1$ we can detect all cosets in $\overline{\mathcal{C}}_n(k,\alpha)$ as in the image of $\mathrm{Tr}_k^P$.
Then we have $c_\alpha=\dim_{\Q}\overline{\mathcal{C}}_n^{\Q}(k,\alpha)$.
Since in this case all necklaces in $\overline{\mathcal{C}}_n(k,\alpha)$ have the preiod $k$
(i.e. primitive), $c_\alpha=\dim_{\Q}\mathcal{L}_n^{\Q}(k,\alpha)$.
Therefore it suffices to consider the case where all positive parts of $\alpha$ are greater than or equal to $2$.
For $5 \leq k \leq 9$, we can calculate $c_\alpha$ and the difference of $r_\alpha=c_\alpha-\dim_{\Q}\mathcal{L}_n^{\Q}(\alpha,k)$ as the following tables.
\[
\begin{array}{cccc}
k & \alpha & c_\alpha & r_\alpha \\
\hline
\hline
5 & (3,2) & 1 & 0 \\
\hline
6 & (4,2) & 2 & 0 \\
& (3,3) & 3 & 0 \\
& (2^3) & 15 & 1 \\
\hline
7 & (5,2) & 3 & 0 \\
 & (4,3) & 5 & 0 \\
 & (3,2^2) & 30 & 0 \\
\end{array} \quad \quad 
\begin{array}{cccc}
k & \alpha & c_\alpha & r_\alpha \\
\hline
\hline
8 & (6,2) & 2 & -1 \\
 & (5,3) & 6 & -1 \\
 & (4,4) & 7 & -1 \\
 & (4,2^2) & 52 & 1 \\
 & (3^2,2) & 69 & -1 \\
 & (2^4) & 316 & 4\\
  & \quad & & 
\end{array}
\quad \quad 
\begin{array}{cccc}
k & \alpha & c_\alpha & r_\alpha \\
\hline
\hline
9 & (7,2) & 4 & 0 \\
 & (6,3) & 9 & 0 \\
 & (5,4) & 14 & 0 \\
 & (5,2^2) & 84 & 0 \\
 & (4,3,2) & 140 & 0 \\
 & (3^3) & 188 & 2 \\
  & (3,2^3) & 840 & 0 
\end{array}
\]
Using above, for $n \geqq 3$, we have
\begin{align*}
&\dim_{\Q}\mathrm{Im}(\overline{\mathrm{Tr}}_{5,\Q}^P)=r_n(5), \quad 
\dim_{\Q}\mathrm{Im}(\overline{\mathrm{Tr}}_{6,\Q}^P)=r_n(6)+\tbinom{n}{3}, \quad
\dim_{\Q}\mathrm{Im}(\overline{\mathrm{Tr}}_{7,\Q}^P)=r_n(7), \\
&\dim_{\Q}\mathrm{Im}(\overline{\mathrm{Tr}}_{8,\Q}^P)=r_n(8)-2n(n-1)-\tbinom{n}{2}+4\tbinom{n}{4} \delta_{n \geq 4}, \\
&\dim_{\Q}\mathrm{Im}(\overline{\mathrm{Tr}}_{9,\Q}^P)=r_n(9)+2\tbinom{n}{3}.
\end{align*}
Notice that $\dim_{\Q}\mathrm{Ker}(\overline{\mathrm{Tr}}_{k,\Q}^P)=n r_n(k)-
\dim_{\Q}\mathrm{Im}(\overline{\mathrm{Tr}}_{k,\Q}^P)$.
Since $\dim_{\Q}\overline{\mathcal{C}}_n^{\Q}(k)=r_n(k)$ for any prime $k$, we remark that $\overline{\mathrm{Tr}}_{k,\Q}^P$ is surjective for $k=5$ and $7$.
But for $k=11$, since we can calcurate that $r_\alpha=-1,-2$ for $\alpha=(8,3),(7,4)$ respectively,
we see $\overline{\mathrm{Tr}}_{11, \Q}^P$ is not surjective.

\vspace{0.5em}

Let us proceed to the observations of the Johnson images.
By fixing $n$, we can directly calcurate of the dimension of the submodule
\[
\mathrm{Im}(\tau_k^P)_1 :=
\langle 
[K_{i_1j_1},K_{i_2j_2}, \ldots ,K_{i_kj_k}] \ | \ 1 \leq i_1,i_2, \ldots i_k, j_1, \ldots ,j_k \leq n
\rangle \subset \mathrm{Im}(\tau_k^P) \subset \mathfrak{p}_n(k)
\]
generated by the degree $1$ part.

In the case of $n=2$, we have $\mathcal{P}_2(k)=\mathcal{A}_2(k)$ for any $k \geq 1$, and see that
the filtration $\mathcal{P}_2(1) \supset \mathcal{P}_2(2) \supset \cdots$ is
the lower central series of the inner automorphism group $\mathrm{Inn}\,F_2$.
Hence $\mathrm{Im}(\tau_k^P) \cong \mathcal{L}_2(k)$ for any $k \geq 1$.
Since $r_2(11)=\mathrm{dim}_{\Q} \mathcal{L}_2(11) < \mathrm{dim}_{\Q} \mathrm{Ker}(\overline{\mathrm{Tr}}_{11,\Q}^P)$ as above,
we have $\mathrm{Im}(\tau_{11,\Q}^P) \subsetneq \mathrm{Ker}(\overline{\mathrm{Tr}}_{11,\Q}^P)$.
In the case of $n=3$, we have the following table on the dimensions of $\mathrm{Im}(\tau_{k,\Q}^P)_1$ and the kernels of $\overline{\mathrm{Tr}}_{k,\Q}^P$.
\[
\begin{array}{c||cccccccc}
k & 1 & 2 & 3 & 4 & 5 & 6 & 7 & 8 \\
\dim_{\Q}\mathrm{Im}(\tau_{k,\Q}^P)_1 & 6 & 6 & 16 & 36 & 96 & 231 & 618 & 1596 \\
\dim_{\Q}\mathrm{Ker}(\overline{\mathrm{Tr}}_{k,\Q}^P) & 6 & 6 & 16 & 36 & 96 & 231 & 624 & 1635
\end{array}
\]
Since $\dim_{\Q}\mathrm{Im}(\tau_{k,\Q}^P)_1=\dim_{\Q}\mathrm{Ker}(\overline{\mathrm{Tr}}_{k,\Q}^P)$ for $1 \leq k \leq 6$, we see
$\mathrm{Im}(\tau_{k,\Q}^P)_1 = \mathrm{Im}(\tau_{k,\Q}^P)=\mathrm{Ker}(\overline{\mathrm{Tr}}_{k,\Q}^P)$ for $1 \leq k \leq 6$.
This shows that the natural homomorphism $\nu : \mathrm{P}\Sigma_3(k)/\mathrm{P}\Sigma_3(k+1) \rightarrow \mathrm{gr}^k(\mathcal{P}_3)$
induced by the inclusion $\mathrm{P}\Sigma_3(k) \hookrightarrow \mathcal{P}_3(k)$ is surjective for $1 \leq k \leq 6$ after tensoring with $\Q$.
But we find a gap at the incusion
$\mathrm{Im}(\tau_{k,\Q}^P)_1 \subsetneq \mathrm{Ker}(\overline{\mathrm{Tr}}_{k,\Q}^P)$ for $k=7$. We can also see that 
this gap appears only the subspace for $(j_1,\ldots ,j_7) \in \mathfrak{S}_3 \cdot (3,2^2)$. 
Moreover, The module $\mathrm{Ker}(\overline{\mathrm{Tr}}_{7,\Q}^P)/\mathrm{Im}(\tau_{7,\Q}^P)_1$ is isomorphic to the regular representation of $\mathfrak{S}_3$.

\vspace{0.5em}

In general, it seems to be difficult to determine whether or not $\mathrm{Im}(\tau_k^P)=\mathrm{Ker}(\overline{\mathrm{Tr}}_k^P)$.
Here we give a partial answer to this problem.
For any composition $\alpha=(\alpha_1,\alpha_2, \ldots ,\alpha_n)$ of $k$,
denote by $\mathfrak{p}_n(k,\alpha)$ the subspace of $\mathfrak{p}_n(k)$ generated by
all $x_i^* \otimes [x_{i_1},x_{i_2}, \ldots ,x_{i_k}, x_i]$ such that $(i_1,i_2, \ldots, i_k) \in p(\alpha)$.
\begin{thm}\label{T-0530}
For any $n \geq 3$, $k \geq 1$ and any
composition $\alpha=(\alpha_1,\alpha_2, \ldots ,\alpha_n)$ of $k$ such that $\alpha_j=1$ for some $1 \leq j \leq n$,
we have $\mathrm{Ker}(\overline{\mathrm{Tr}}_k^P|_{\mathfrak{p}_n(k,\alpha)}) \subset \mathrm{Im}(\tau_k^P)$.
\end{thm}
\textit{Proof.}
For any $X \in \mathrm{Ker}(\overline{\mathrm{Tr}}_k^P|_{\mathfrak{p}_n(k,\alpha)})$, set
\[ X = \sum a_{i, i_1, \ldots, i_k} x_i^* \otimes [x_{i_1},x_{i_2}, \ldots, x_{i_k}, x_i] \]
for some $a_{i, i_1, \ldots, i_k} \in \Z$.

\vspace{0.5em}

\noindent
{\bf Step 1}. We show that $x_i^* \otimes [x_{i_1}, \ldots, x_{i_k}, x_i] \in \mathrm{Im}(\tau_k^P)$ if $i_1, \ldots, i_k \neq i$
by the induction on $k \geq 1$. For $k =1$, it is trivial since $\tau_1^P$ is surjective. Assume $k \geq 2$.
By the inductive hypothesis, there exist some $\sigma \in \mathcal{P}_n(k-1)$ and $\tau \in \mathcal{P}_n(1)$ such that
\[ \tau_{k-1}^P(\sigma) = x_i^* \otimes [x_{i_1},x_{i_2}, \ldots, x_{i_{k-1}}, x_i],
   \hspace{1em} \tau_{1}^P(\tau) = x_i^* \otimes [x_{i_k}, x_i]. \]
Then we have $\tau_k^P([\sigma, \tau])=x_i^* \otimes [x_{i_1},x_{i_2}, \ldots, x_{i_k}, x_i]$, and hence the induction proceeds.
Since $\mathrm{Im}(\tau_k^P) \subset \mathrm{Ker}(\overline{\mathrm{Tr}}_k^P)$, 
this shows that we may assume that $X$ is a linear combination of $x_i^* \otimes [x_{i_1}, \ldots, x_{i_k}, x_i]$
such that there exists some $1 \leq l \leq k$ such that $i_l=i$.
Furthermore, if we use the Jacobi identity
repeatedly, we have
\[\begin{split}
   x_i^* & \otimes [x_{i_1}, \ldots x_{i_{l-1}}, x_i, x_{i_{l+1}}, \ldots, x_{i_k}, x_i] \\
   & =- x_i^* \otimes [x_i, [x_{i_1}, \ldots x_{i_{l-1}}], x_{i_{l+1}}, \ldots, x_{i_k}, x_i], \\
   & =- x_i^* \otimes [x_i, [x_{i_1}, \ldots x_{i_{l-2}}], x_{i_{l-1}}, x_{i_{l+1}}, \ldots, x_{i_k}, x_i]
      + x_i^* \otimes [x_i, x_{i_{l-1}}, [x_{i_1}, \ldots x_{i_{l-2}}], x_{i_{l+1}}, \ldots, x_{i_k}, x_i] \\
   & = \cdots.
  \end{split}\]
Thus we may assume that $X$ is a linear combination of elements of type
$x_i^* \otimes [x_i, x_{i_2}, \ldots, x_{i_k}, x_i]$.

\vspace{0.5em}

\noindent
{\bf Step 2}. We consider elements $x_i^* \otimes [x_i, x_{i_2}, \ldots, x_{i_k}, x_i]$ such that $i_p=i$ for some $2 \leq p \leq k$.
By the hypothesis of the proposition for $\alpha$, there exists some $2 \leq l \leq k$ such that $i_l$ appears only one time in $i, i_2, \ldots, i_k$.
We see $i_l \neq i$. We show that $x_i^* \otimes [x_i, x_{i_2}, \ldots, x_{i_k}, x_i]$ is a linear combination of elements of type
\[ x_{i_l}^* \otimes [x_{i_l}, x_{j_2}, \ldots, x_{j_k}, x_{i_l}] \hspace{1em} (j_2, \ldots, j_k \neq i_l) \]
modulo $\mathrm{Im}(\tau_k^P)$ by the induction on $k \geq 2$.
For $k=2$, we consider $x_i^* \otimes [x_i, x_{i_2}, x_i]$ for $i_2 \neq i$. We see that the claim hold in this case by
\[ \tau_2^P([K_{i_2 i}, K_{i i_2}]) = x_{i_2}^* \otimes [x_{i_2}, x_i, x_{i_2}] - x_i^* \otimes [x_i, x_{i_2}, x_i]. \]
Assume $k \geq 3$. First we consider the case where $i_{l+1}, \ldots, i_k \neq i$.
From the argument is Step 1, we see
\[\begin{split}
   \mathrm{Im}(\tau_k^P) \ni & [x_i^* \otimes [x_{i_l}, x_{i_{l+1}}, \ldots, x_{i_k}, x_i], x_{i_l}^* \otimes [x_i, x_{i_2}, \ldots, x_{i_{l-1}}, x_{i_l}]] \\
   & = x_i^* \otimes [x_i, x_{i_2}, \ldots, x_{i_k}, x_i] - x_{i_l}^* \otimes [x_{i_l}, x_{i_{l+1}}, \ldots, x_{i_k}, x_i, x_{i_2}, \ldots, x_{i_l}] \\
   & \hspace{0.5em} - \sum_{\substack{2 \leq m \leq l\\[1pt] i_m=i}} x_{i_l}^* \otimes [x_{i}, x_{i_2}, \ldots, x_{i_{m-1}}, [x_{i_l}, \ldots, x_{i_k}, x_i],
                         x_{i_{m+1}}, \ldots, x_{i_l}],
  \end{split}\]
and see that the claim holds by using the Jacobi identity.
Next we consider the case where $i_m=i$ for some $l+1 \leq m \leq k$. By the inductive hypothesis,
an element
\[ x_i^* \otimes [x_{i_l}, x_{i_{l+1}}, \ldots, x_{i_k}, x_i] = - x_i^* \otimes [x_i, [x_{i_l}, \ldots, x_{i_{m-1}}], x_{i_{m+1}}, \ldots, x_{i_k}, x_i] \]
can be written as a linear combination of elements of type
$x_{i_l}^* \otimes [x_{i_l}, x_{j_2}, \ldots, x_{j_{k-l+1}}, x_{i_l}]$ for $j_2, \ldots, j_{k-l+1} \neq i_l$ modulo $\mathrm{Im}(\tau_{k-l+1}^P)$.
On the other hand, we see
\[\begin{split}
  [x_{i_l}^* \otimes [x_{i_l}, & x_{j_2}, \ldots, x_{j_{k-l+1}}, x_{i_l}], x_{i_l}^* \otimes [x_i, x_{i_2}, \ldots, x_{i_{l-1}}, x_{i_l}]] \\
   & = x_{i_l}^* \otimes [x_i, x_{i_2}, \ldots, x_{i_{l-1}}, x_{i_l}, x_{j_2}, \ldots, x_{j_{k-l+1}}, x_{i_l}] \\
   & \hspace{1em} + x_{i_l}^* \otimes [x_{i_l}, x_{j_2}, \ldots, x_{j_{k-l+1}}, [x_i, x_{i_2}, \ldots, x_{i_{l-1}}, x_{i_l}]] \\
   & \hspace{1em} - x_{i_l}^* \otimes [x_i, x_{i_2}, \ldots, x_{i_{l-1}}, [x_{i_l}, x_{j_2}, \ldots, x_{j_{k-l+1}}, x_{i_l}]] \\
   & = x_{i_l}^* \otimes [x_i, x_{i_2}, \ldots, x_{i_{l-1}}, x_{i_l}, x_{j_2}, \ldots, x_{j_{k-l+1}}, x_{i_l}] \\
   & \hspace{1em} + x_{i_l}^* \otimes [x_{i_l}, x_{j_2}, \ldots, x_{j_{k-l+1}}, [x_i, x_{i_2}, \ldots, x_{i_{l-1}}], x_{i_l}].
  \end{split}\]
Therefore by the same arguments as the previous case, we see that the claim holds by using the Jacobi identity.
Thus the induction proceeds.
we may assume that $X$ is a linear combination of elements of type
$x_i^* \otimes [x_i, x_{i_2}, \ldots, x_{i_k}, x_i]$ for $i_2, \ldots, i_k \neq i$.

\vspace{0.5em}

\noindent
{\bf Step 3}. Consider an element $x_i^* \otimes [x_i, x_{i_2}, \ldots, x_{i_k}, x_i]$ for $i_2, \ldots, i_k \neq i$.
If there exists another $2 \leq l \leq k$ such that $i_l$ appears only one time in $i, i_2, \ldots, i_k$,
from the argument in Step 1, we have some $\sigma \in \mathcal{P}_n(k-l+1)$ and $\tau \in \mathcal{P}_n(l)$ such that
\[ \tau_{k-l+1}^P(\sigma) = x_i^* \otimes [x_{i_l},x_{i_{l+1}}, \ldots, x_{i_k}, x_i],
   \hspace{1em} \tau_{l-1}^P(\tau) = x_{i_l}^* \otimes [x_i, x_{i_2}, \ldots, x_{i_l}], \]
and
\[ \tau_k^P([\sigma,\tau]) = x_i^* \otimes [x_i, x_{i_2}, \ldots, x_{i_k}, x_i] 
     - x_{i_l}^* \otimes [x_{i_l}, x_{i_{l+1}}, \ldots, x_{i_k}, x_i, x_{i_2}, \ldots, x_{i_l}]. \]
This shows that we may assume that $X$ is a linear combination of elements of type
$x_i^* \otimes [x_i, x_{i_2}, \ldots, x_{i_k}, x_i]$ for $i_2, \ldots, i_k \neq i$ such that $i$ is the minimum element in $i, i_2, \ldots, i_k$,
which appears only one time in $i, i_2, \ldots, i_k$. Set
\[ X = {\sum}' a_{i, i_2, \ldots, i_k} x_i^* \otimes [x_i, x_{i_2}, \ldots, x_{i_k}, x_i]  \]
where the sum runs over all $(i, i_2, \ldots, i_k)$ such that $i_2, \ldots, i_k \neq i$ and $i$ is the minimum element in $i, i_2, \ldots, i_k$,
which appears only one time in $i, i_2, \ldots, i_k$.
Then we have
\[ \overline{\mathrm{Tr}}_k^P(X) = {\sum}' a_{i, i_2, \ldots, i_k} x_i x_{i_2} \cdots x_{i_k} \in \overline{\mathcal{C}}_n(k). \]
Since elements $x_i x_{i_2} \cdots x_{i_k}$ are linearly independent in $\overline{\mathcal{C}}_n(k)$, we obtain
$a_{i, i_2, \ldots, i_k}=0$ for all $(i, i_2, \ldots, i_k)$.
This completes the proof of Theorem {\rmfamily \ref{T-0530}}.
$\square$

\section{The Andreadakis conjecture for $\mathfrak{p}_n(4)$}

In this section, as an application of Proposition {\rmfamily \ref{P-2}},
we show that $\mathcal{P}_n(4)=\mathrm{P}\Sigma_n(4)$. 
First, we give the strategy to show it, which is the same way as the proof of $\mathcal{P}_n(3)=\mathrm{P}\Sigma_n(3)$ given in \cite{S15}.
Since the natural map
\[ \nu : \mathrm{P}\Sigma_n(3)/\mathrm{P}\Sigma_n(4) \rightarrow \mathrm{gr}^3(\mathcal{P}_n) \]
is surjective, and $\mathrm{P}\Sigma_n(4) \subset \mathrm{P}\Sigma_n(4)$, if the injectivity of $\nu$ is proved then $\nu$ is an isomorphism, and hence
we obtain $\mathrm{P}\Sigma_n(4) = \mathcal{P}_n(4)$.
Since $\mathrm{gr}^3(\mathcal{P}_n)$ is a free abelian group of rank $n(n-1)^2(n+1)/3$, if we show that 
\begin{center}
\lq\lq $\mathrm{P}\Sigma_n(3)/\mathrm{P}\Sigma_n(4)$ is generated by $n(n-1)^2(n+1)/3$ elements as a $\Z$-module"
\end{center}
\noindent
then there exists a section
$\mathrm{gr}^3(\mathcal{P}_n) \rightarrow \mathrm{P}\Sigma_n(3)/\mathrm{P}\Sigma_n(4)$ which is the inverse of $\nu$, and hence
$\nu$ is injective. Thus we show the following.
\begin{thm}\label{T-main-3}
$\mathrm{P}\Sigma_n(3)/\mathrm{P}\Sigma_n(4)$ is generated by \\
\hspace{2em} {\bf{(E1)}}: $[K_{ij}, K_{il}, K_{im}]$ for any distinct $1 \leq i, j, l, m \leq n$ such that $j>l<m$, \\
\hspace{2em} {\bf{(E2)}}: $[K_{ij}, K_{il}, K_{ij}]$ for any distinct $1 \leq i, j, l \leq n$, \\
\hspace{2em} {\bf{(E3)}}: $[K_{ij}, K_{il}, K_{ji}]$ for any distinct $1 \leq i, j, l \leq n$ such that not $i>j, k$, \\
\hspace{2em} {\bf{(E4)}}: $[K_{ij}, K_{ji}, K_{ij}]$ for any distinct $1 \leq i, j \leq n$
as a $\Z$-module. The number of the above elements is $n(n-1)^2(n+1)/3$.
\end{thm}
\textit{Proof.}
First, we count the number of the elements.
The number of elements of {\bf{(E1)}} and {\bf{(E2)}} is equal to the total number of choices of $i \in \{ 1, \ldots, n\}$ times the total number of choices
of $j, l, m \in \{ 1 \ldots, i-1, i+1, \ldots, n\}$ such that $j>l \leq m$. It is $n \times 2 \binom{(n-1)+1}{3}$.
The number of elements of {\bf{(E3)}} is $\binom{n}{3} \times 4$. The number of elements of {\bf{(E4)}} is $n(n-1)$.
Thus the total number of the elements is equal to $n(n-1)^2(n+1)/3$.

\vspace{0.5em}

From (\ref{eq-basis-tau_2}), we see that
$\mathrm{P}\Sigma_n(2)/\mathrm{P}\Sigma_n(3)=\mathrm{gr}^2(\mathcal{P}_n)$ is generated by
\[\begin{split}
    \{ [K_{ij}, K_{il}] \,|\, 1 \leq l < j \leq n, \hspace{0.5em} i \neq j, l \} \cup \{ [K_{ij}, K_{ji}] \,|\, 1 \leq i \neq j \leq n \}.
  \end{split} \]
Hence $\mathrm{P}\Sigma_n(3)/\mathrm{P}\Sigma_n(4)$ is generated by
\[\begin{cases}
   [K_{ij}, K_{il}, K_{pq}] \hspace{1em} & i \neq l<j \neq i, \\
   [K_{ij}, K_{ji}, K_{pq}] & i \neq j
  \end{cases}\]
for any $p \neq q$. We reduce these generators.

\vspace{0.5em}

\noindent
(I) The generators type of $X:=[K_{ij}, K_{il}, K_{pq}]$.

Set $N:=\sharp \{ i, j, l, p, q \}$. If $N=5$, then $X=0$ since as an automorphism of $F_n$, $[K_{ij}, K_{il}, K_{pq}]=1$ in $\mathrm{Aut}\,F_n$.
Consider the case where $N=4$. In this case, $p \in \{ i, j, l \}$ or $q \in \{i, j, l \}$.
If $p=i$ and $l<q$, then $X$ is of type {\bf{(E1)}}. If $p=i$ and $l>q$, by using the Jacobi identity, we see
\[ [K_{ij}, K_{il}, K_{iq}] = -[K_{il}, K_{iq}, K_{ij}] +[K_{ij}, K_{iq}, K_{il}]. \]
If $p=j$, by the Jacobi idenity and $[K_{ij}, K_{jq}]=-[K_{ij}, K_{iq}]$ from the relation {\bf{(P3)}}, we have
\[ [K_{ij}, K_{il}, K_{jq}] = -[K_{il}, K_{jq}, K_{ij}] -[K_{jq}, K_{ij}, K_{il}] = [K_{ij}, K_{jq}, K_{il}] = -[K_{ij}, K_{iq}, K_{il}]. \]
This is the same case as $p=i$.
Hence, in each of the above cases, $X$ can be written as a sum of elements type of {\bf{(E1)}}.
So, we can remove $X$ from the generating set in this case.
Similarly, we can reduce the generators in other cases. For simplicity, we give a list as follows.
\vspace{0.5em}
\begin{center}
{\renewcommand{\arraystretch}{1.1}
\begin{tabular}{|c|l|l|l|} \hline
  $N$  &               & $X$ is equal to & a linear combination of  \\ \hline
  $5$  &               & $0$ &   \\ \hline
  $4$  & $p=i$         & $[K_{ij}, K_{il}, K_{iq}]$ & {\bf{(E1)}} \\ \hline 
       & $p=j$         & $[K_{iq}, K_{ij}, K_{il}]$ & {\bf{(E1)}} \\ \hline
       & $p=l$         & $[K_{il}, K_{iq}, K_{ij}]$ & {\bf{(E1)}} \\ \hline
       & $q=i$         & $[K_{pl}, K_{pj}, K_{pi}]$ & {\bf{(E1)}} \\ \hline
       & $q=j$         & $0$ &  \\ \hline
       & $q=l$         & $0$ &  \\ \hline
\end{tabular}}
\end{center}
\vspace{0.5em}
\begin{center}
{\renewcommand{\arraystretch}{1.1}
\begin{tabular}{|c|l|l|l|} \hline
  $N$  &                 & $X$ is equal to & a linear combination of  \\ \hline
  $3$  & $(p,q)=(i,j)$   & $[K_{ij}, K_{il}, K_{ij}]$  & {\bf{(E2)}} \\ \hline 
       & $(p,q)=(i,l)$   & $-[K_{il}, K_{ij}, K_{il}]$ & {\bf{(E2)}} \\ \hline
       & $(p,q)=(j,i)$   & $[K_{ij}, K_{il}, K_{ji}]$  & {\bf{(E3)}} \\ \hline
       & $(p,q)=(j,l)$   & $[K_{il}, K_{ij}, K_{il}]$ & {\bf{(E2)}} \\ \hline
       & $(p,q)=(l,i)$   & $-[K_{il}, K_{ij}, K_{li}]$ & {\bf{(E3)'}} \\ \hline
       & $(p,q)=(l,j)$   & $-[K_{ij}, K_{il}, K_{ij}]$ & {\bf{(E2)}} \\ \hline
\end{tabular}}
\end{center}
\vspace{0.5em}
\noindent
Here {\bf{(E3)'}} means
\begin{center}
 {\bf{(E3)'}}: $[K_{ij}, K_{il}, K_{ji}]$ for any distinct $1 \leq i, j, l \leq n$.
\end{center}
\noindent
We show that any element type of {\bf{(E3)'}} can be written as that of {\bf{(E3)}} later.

\vspace{0.5em}

\noindent
(II) The generators type of $Y:=[K_{ij}, K_{ji}, K_{pq}]$.

Set $N':=\sharp \{ i, j, p, q \}$. By arguments similar to those in (I), we have the following.
\vspace{0.5em}
\begin{center}
{\renewcommand{\arraystretch}{1.1}
\begin{tabular}{|c|l|l|l|} \hline
  $N'$ &               & $Y$ is equal to & a linear combination of  \\ \hline
  $4$  &               & $0$ &   \\ \hline
  $3$  & $p=i$         & $[K_{ji}, K_{jq}, K_{ij}]+ [K_{ij}, K_{iq}, K_{ji}]$ & {\bf{(E3)'}} \\ \hline 
       & $p=j$         & $-[K_{ij}, K_{iq}, K_{ji}]- [K_{ji}, K_{jq}, K_{ij}]$ & {\bf{(E3)'}} \\ \hline
       & $q=i$         & $-[K_{pi}, K_{pj}, K_{pi}]$ & {\bf{(E2)}} \\ \hline
       & $q=j$         & $[K_{pj}, K_{pi}, K_{pj}]$ & {\bf{(E2)}} \\ \hline
  $2$  & $(p,q)=(i,j)$ & $[K_{ij}, K_{ji}, K_{ij}]$ & {\bf{(E4)}} \\ \hline
       & $(p,q)=(j,i)$ & $-[K_{ji}, K_{ij}, K_{ji}]$ & {\bf{(E4)}} \\ \hline
\end{tabular}}
\end{center}
\vspace{0.5em}

Finally we consider {\bf{(E3)'}}. For any distinct $1 \leq i, j, l \leq n$, we have 
\[\begin{split}
   0 & \stackrel{\text{{\bf{(P3)}}}}{=} [K_{ij}K_{lj}, K_{il}, K_{ji}K_{li}] \\
     & = [K_{ij}, K_{il}, K_{ji}K_{li}] + \underline{[K_{lj}, K_{il}, K_{ji}K_{li}]} \\
     & \stackrel{\text{J, {\bf{(P3)}}}}{=} [K_{ij}, K_{il}, K_{ji}K_{li}] - \underline{[K_{il}, K_{ji}K_{li}, K_{lj}]} \\
     & = [K_{ij}, K_{il}, K_{ji}] + \uwave{[K_{ij}, K_{il}, K_{li}]} - \uuline{[K_{il}, K_{ji}, K_{lj}]} - [K_{il}, K_{li}, K_{lj}] \\
     & \stackrel{\text{{\bf{(P3)}}}}{=} [K_{ij}, K_{il}, K_{ji}] + \uwave{[K_{il}, K_{lj}, K_{li}]} - \uuline{[K_{ji}, K_{jl}, K_{lj}]} - [K_{il}, K_{li}, K_{lj}] \\
     & \stackrel{\text{J}}{=} [K_{ij}, K_{il}, K_{ji}] + [K_{jl}, K_{ji}, K_{lj}] + [K_{li}, K_{lj}, K_{il}]
 \end{split}\]
where $\stackrel{\text{J}}{=}$ and $\stackrel{\text{{\bf{(P3)}}}}{=}$ mean the equlities induced from the Jacobi identity and
the relation {\bf{(P3)}} respectively. Hence for any $1 \leq l<j<i \leq n$,
we can remove $[K_{ij}, K_{il}, K_{ji}]$ and $[K_{il}, K_{ij}, K_{li}]$ from the generating set.

\vspace{0.5em}

Therefore we obtain the required result. This completes the proof of Theorem {\rmfamily \ref{T-main-3}}. $\square$

\vspace{0.5em}

Hence we obtain the following.
\begin{cor}\label{C-6}
For $n \geq 3$, we have $\mathcal{P}_n(4)=\mathrm{P}\Sigma_n(4)$.
\end{cor}

\section{Remarks on $H_2(\mathrm{P}\Sigma_n, \Z)$}

In this section, as an applications of our results above,
we give some remarks on the second homology and the second cohomology groups of $\mathrm{P}\Sigma_n$.
All of the results in this subsection were alreday obtained by Brownstein and Lee \cite{BrL}, and
Jensen, McCammond and Meiyer \cite{JMM}.
The method is the same as that in the study of the integral second homology group of the lower-triangular IA-automorphism group of $F_n$ in \cite{S05}.

\vspace{0.5em}

We prepare some notation. Let $F$ be the free group of rank $n(n-1)$ with a basis $K_{ij}$ for all $1 \leq i \neq j \leq n$.
Let $\pi : F \rightarrow \mathrm{P}\Sigma_n$ be the standard surjective homomorphism, and set $R:=\mathrm{Ker}(\pi)$.
Then we have the group extension
\begin{equation}\label{eq-ext-pres}
1 \rightarrow R \rightarrow F \xrightarrow{\pi} \mathrm{P}\Sigma_n \rightarrow 1.
\end{equation}
Since the abelianization of $\mathrm{P}\Sigma_n$ is a free abelian group generated by
(the coset classes of) all $K_{ij}$s, we see that $\pi$ induces an isomorphism
\[ \pi_* : H_1(F,\Z) \rightarrow H_1(\mathrm{P}\Sigma_n,\Z). \]
From the homological five-term exact sequence of (\ref{eq-ext-pres}),
we have
\[\begin{split}
   0=H_2(F,\Z) \rightarrow H_2(\mathrm{P}\Sigma_n,\Z) \rightarrow H_1(R,\Z)_F 
    \rightarrow H_1(F,\Z) \xrightarrow{\pi_*} H_1(\mathrm{P}\Sigma_n,\Z) \rightarrow 0,
  \end{split}\]
and hence
\[ H_2(\mathrm{P}\Sigma_n,\Z) \cong H_1(R,\Z)_F. \]

\vspace{0.5em}

Let $F=\Gamma_F(1) \supset \Gamma_F(2) \supset \cdots$ be the lower central series of $F$, and set
$\mathcal{L}_F(k) := \Gamma_F(k)/\Gamma_F(k+1)$ for each $k \geq 1$.
Let $\{ R_k \}_{k \geq 1}$ be the descending filtration of $R$ defined by $R_k := R \cap \Gamma_F(k)$ for each $k \geq 1$. 
We have $R_k=R$ for $k=1$ and $2$ since $R \subset [F,F]$.
For each $k \geq 1$, let
$\pi_k : \mathcal{L}_F(k) \rightarrow \mathrm{P}\Sigma_n(k)/\mathrm{P}\Sigma_n(k+1)$
be the homomorphism induced from $\pi : F \rightarrow \mathrm{P}\Sigma_n$.
By observing $R_k/R_{k+1} \cong (R_k \, \Gamma_F(k+1))/\Gamma_F(k+1)$, we obtain an exact sequence
\begin{equation}\label{ZZ}
0 \rightarrow R_k/R_{k+1} \rightarrow \mathcal{L}_F(k) \xrightarrow{\pi_k} \mathrm{P}\Sigma_n(k)/\mathrm{P}\Sigma_n(k+1) \rightarrow 0.
\end{equation}
For each $k \geq 2$.
The natural projection $R \rightarrow R/R_{k+1}$ induces the surjective homomorphism
$\psi_k : H_1(R,\Z) \rightarrow H_1(R/R_{k+1},\Z)$.
By considering the right action of $F$ on $R$, defined by
\[ r \cdot x := x^{-1} r x, \hspace{1em} r \in R, \,\,\, x \in F, \]
we see that $\psi_k$ is $F$-equivariant.
Hence it induces the surjective homomorphism
\[ H_1(R,\Z)_{F} \rightarrow H_1(R/R_{k+1},\Z)_{F}, \]
which is also denoted by $\psi_k$.
For $k=2$, $H_1(R/R_{3},\Z)_{F} = R/R_{3}$ since $F$ acts on $R/R_{3}$ trivially.

\begin{pro}\label{T-main-4}
For $n \geq 3$, $H_2(\mathrm{P}\Sigma_n,\Z) \cong \Z^{\frac{1}{2}n^2(n-1)(n-2)}$.
\end{pro}
\textit{Proof.}
Observe $\mathrm{P}\Sigma_n(2)/\mathrm{P}\Sigma_n(3)=\mathrm{gr}^2(\mathcal{P}_n)$ is a free abelian group of $n(n-1)^2/2$ as mentioned above,
by considering (\ref{ZZ}) for $k=2$, we see
\[\begin{split}
   \mathrm{rank}_{\Z}(R/R_{3}) & = \mathrm{rank}_{\Z}(\mathcal{L}_F(2)) - \mathrm{rank}_{\Z}(\mathrm{gr}^2(\mathcal{P}_n)) \\
    & = \frac{1}{2}n(n-1)(n(n-1)-1) - \frac{1}{2}n(n-1)^2 \\
    & = \frac{1}{2}n^2(n-1)(n-2).
   \end{split}\]
This means $H_2(\mathrm{P}\Sigma_n,\Z)$ contains a free abelian group of rank $n^2(n-1)(n-2)/2$.
On the other hand, $H_1(R,\Z)_F \cong R/[F,R]$ is generated by {\bf{(P1)}}, {\bf{(P2)}} and {\bf{(P3)}} as a $\Z$-module.
Since the total number of elements of {\bf{(P1)}}, {\bf{(P2)}} and {\bf{(P3)}} is $n^2(n-1)(n-2)/2$, we see that
$\psi_2 : R/[F,R] = H_1(R,\Z)_F \rightarrow H_1(R/R_3,\Z)_{F}=R/R_3$ is an isomorphism, and obtain the required result.
$\square$

\vspace{0.5em}

By the universal coefficient theorem, we have
\[ H^2(\mathrm{P}\Sigma_n,\Z) \cong \Z^{\frac{1}{2}n^2(n-1)(n-2)}. \]
\begin{pro}\label{T-main-5}
For $n \geq 3$, the cup product map $\cup : \wedge^2 H^1(\mathrm{P}\Sigma_n, \Z) \rightarrow H^2(\mathrm{P}\Sigma_n, \Z)$ is surjective.
\end{pro}
\textit{Proof.}
By observing the cohomological five-term exact sequence of (\ref{eq-ext-pres}), we can see that
$H^2(\mathrm{P}\Sigma_n,\Z) \cong H^1(R,\Z)^F=H^1(R/[F,R],\Z)$. From the argument in the proof of Proposition {\rmfamily \ref{T-main-4}},
since the natural map $R/[F,R] \rightarrow R/R_3$ is an isomorphism, $H^2(\mathrm{P}\Sigma_n,\Z) \cong H^1(R/R_3,\Z)$.
Since $H^1(R/R_3,\Z)$ coincides with the image of the cup product map, we obtain the required result.
(For details, see Lemma 4.1 in \cite{S06}.) $\square$

\section{Twisted first cohomology groups of $\mathrm{BP}_n$}

From the group extension
\[ 1 \rightarrow \mathrm{P}\Sigma_n \rightarrow \mathrm{BP}_n \rightarrow \mathfrak{S}_n \rightarrow 1, \]
the symmetric group $\mathfrak{S}_n$ acts on the integral homology groups $H_p(\mathrm{P}\Sigma_n, \Z)$ for $p \geq 1$.
We consider $H_p(\mathrm{P}\Sigma_n, \Z)$ as a $\mathrm{BP}_n$-module through the homomorphism $\mathrm{BP}_n \rightarrow \mathfrak{S}_n$.

Note that $\mathfrak{p}_n^{\Q}(1)\cong H_{\Q} \otimes \mathcal{L}_n^{\Q}(1) \cong ((n) \oplus (n-1,1)) \otimes (n-1,1)$ as an $\mathfrak{S}_n$-module.
We can decompose the tensor product representation $(n-1,1) \otimes (n-1,1) \cong (n) \oplus (n-1,1) \oplus (n-2,2) \oplus (n-2,1^2)$. (See the example in \cite[1.8]{Tok}).
Therefore, as an $\mathfrak{S}_n$-equivariant $\Q$-vector space we have
\[ H_1(\mathrm{P}\Sigma_n, \Q) =\mathfrak{p}_n^{\Q}(1) \cong (n) \oplus 2(n-1,1) \oplus (n-2,1^2) \oplus (n-2,2) \]
for $n \geq 4$. In this section, we consider the first cohomology groups of $\mathrm{BP}_n$ with coefficients in $(n)$ and $(n-1,1)$.

\vspace{0.5em}

First, the irreducible $\mathfrak{S}_n$-module $(n)$. This is isomorphic to the trivial module $\Q$.
From Theorem {\rmfamily \ref{T-FRR}} by Fenn-Rim\'{a}nyi-Rourke, we see
\[ H_1(\mathrm{BP}_n, \Z) = \Z \oplus \Z/2\Z. \]
The free part is generated by (the coset class of) $\sigma_1$, and the torsion part is generated by (the coset class of) $s_1$.
Thus, by the universal coefficient theorem, we obtain $H^1(\mathrm{BP}_n, \Z) \cong \Z$, and $H^1(\mathrm{BP}_n, \Q) \cong \Q$.

\vspace{0.5em}

Next, we consider the standard $\mathfrak{S}_n$-representation $(n-1,1)$. Recall the standard $\mathrm{GL}(n,\Z)$-representation $H$.
Let $\bm{e}_1, \ldots, \bm{e}_n \in H$ be the standard basis of $H$. For any element $x=x_1\bm{e}_1 + \cdots + x_n \bm{e}_n \in H$,
we write $x=(x_1, \ldots, x_n)$ for simplicity.
As an $\mathfrak{S}_n$-representation, the module $H$ is decomposed into the direct sum of the two irreducible representations
\[\begin{split}
   V &:=\{ (x_1, \ldots, x_n) \in H \,|\, x_1+ \cdots + x_n=0 \}, \\
   W &:=\{ (x, \ldots, x) \in H \,|\, x \in \Z \}.
  \end{split}\] 
The representation $W$ is the trivial representation, and $V \otimes_{\Z} \Q \cong (n-1,1)$.
We compute the first cohomology group of $\mathrm{BP}_n$ with coefficients in $V$. Here we consider $\bm{e}_1-\bm{e}_n, \ldots, \bm{e}_{n-1}-\bm{e}_n$
as a basis of $V$, and write any element $v=v_1(\bm{e}_1-\bm{e}_n) + \cdots + v_{n-1}(\bm{e}_{n-1}-\bm{e}_n) \in V$ as
$v=(v_1, \ldots, v_{n-1})$ for simplicity.

\vspace{0.5em}

Here we give the outline of the computation. 
Let $Z^1(\mathrm{BP}_n, V)$ be the $\Z$-module of all crossed homomorphisms from $\mathrm{BP}_n$ to $V$, and 
$B^1(\mathrm{BP}_n, V)$ be the $\Z$-module of all principal crossed homomorphisms from $\mathrm{BP}_n$ to $V$.
Then we have $H^1(\mathrm{BP}_n, V) \cong Z^1(\mathrm{BP}_n, V)/B^1(\mathrm{BP}_n, V)$.
We compute $Z^1(\mathrm{BP}_n, V)$ by using the presentation of $\mathrm{BP}_n$ in Theorem {\rmfamily \ref{T-FRR}}.
Let $F$ be the free group with a basis $\sigma_i$ and $s_i$ for $1 \leq i \leq n-1$,
and $\varphi : F \rightarrow \mathrm{BP}_n$ the standard projection.
Then the kernel $R$ of $\varphi$ is the normal closure of all of the relators coming from the relations in Theorem {\rmfamily \ref{T-FRR}}.
Considering the five-term exact sequence of the Lyndon-Hochshild-Serre spectral sequence of the group extension
\[ 1 \rightarrow R \rightarrow F \rightarrow \mathrm{BP}_n \rightarrow 1, \]
we obtain the exact sequence
\[ 0 \rightarrow H^1(\mathrm{BP}_n, V) \rightarrow H^1(F, V) \rightarrow H^1(R, V)^F. \]
By observing this sequence at the cocycle level, we also obtain the exact sequence
\begin{equation}\label{eq-cross}
 0 \rightarrow Z^1(\mathrm{BP}_n,V) \rightarrow Z^1(F,V)
       \xrightarrow{\iota^*} Z^1(R,V)
\end{equation}
where $\iota^*$ is the map induced from the inclusion $\iota : R \hookrightarrow F$.
Thus we consider $Z^1(\mathrm{BP}_n,V)$ as a subgroup of $Z^1(F,V)$
which are killed by $\iota^*$. In other words, $Z^1(\mathrm{BP}_n,V)$ consists of all elements of $Z^1(F,V)$ which preserve all of the relations in Theorem {\rmfamily \ref{T-FRR}}.
By the universality of the free group $F$, it is easily seen that a crossed homomorphism $f: F \rightarrow V$ is uniquely determined by the images of the basis
$\sigma_i$ and $s_i$ for $1 \leq i \leq n-1$, and hence $Z^1(F,V)$ is isomorphic to $V^{\oplus 2(n-1)}$.

\begin{thm}\label{T-main-6}
For $n \geq 3$, $H^1(\mathrm{BP}_n, V) \cong \Z^{\oplus 2} \oplus \Z/4\Z$.
\end{thm}
\textit{Proof.}
First, we determine $B^1(\mathrm{BP}_n,V)$. For any $v=(v_1, \ldots, v_{n-1}) \in V$, consider the principal crossed homomorphism $f_v : \mathrm{BP}_n \rightarrow V$
associated to $v$. The images of $\sigma_i$ and $s_i$ for $1 \leq i \leq n-1$ by $f_v$ are given as follows.
\[\begin{split}
    f_v(\sigma_i) & = \begin{cases}
                        (0, \ldots, 0, v_{i+1}-v_i, v_i-v_{i+1}, 0, \ldots, 0) \hspace{1em} & i \neq n-1, \\
                        (0, \ldots, 0, -(v_1+ \cdots+ v_{n-2} + 2v_{n-1})) & i=n-1,
                      \end{cases} \\
    f_v(s_i) & = f_v(\sigma_i).
  \end{split}\]
Here $v_{i+1}-v_i$ and $v_i-v_{i+1}$ are in the $i$-th and the $(i+1)$-st component of $f_v$.

\vspace{0.5em}

Next, we consider $Z^1(\mathrm{BP}_n,V)$. Let $f : F \rightarrow V$ be a crossed homomorphism
such that
\[ f(\sigma_i)=(a_1(i), \ldots, a_{n-1}(i)), \hspace{1em} f(s_i)=(b_1(i), \ldots, b_{n-1}(i)) \]
for any $1 \leq i \leq n-1$.
It suffices to find a condition among $a_j(i)$s and $b_j(i)$s such that $f$ preserves all relations in Theorem {\rmfamily \ref{T-FRR}}.

\vspace{0.5em}

\noindent
{\bf{(B1)}}: From the ralation {\bf{(B1)}} for any $1 \leq i \leq n-2$,
we have $f(\sigma_i \sigma_{i+1} \sigma_i) = f(\sigma_{i+1} \sigma_{i} \sigma_{i+1})$. It is equivalent to
\[ (1-\sigma_i\sigma_{i+1}-\sigma_{i+1})f(\sigma_i)=(1-\sigma_{i+1}\sigma_i-\sigma_i)f(\sigma_{i+1}).  \]
For $i \leq n-3$, we have
{\small
\[\begin{split}
   (\text{LHS})&=(a_1(i), \ldots, a_{i-1}(i), a_{i+2}(i), a_{i+1}(i)-a_{i+2}(i)+a_i(i), a_{i+2}(i), \ldots, a_{n-1}(i)), \\
   (\text{RHS})&=(a_1(i+1), \ldots, a_{i-1}(i+1), a_i(i+1), a_{i+1}(i+1)-a_{i}(i+1)+a_{i+2}(i+1), a_i(i+1), \\
               & \hspace{3em} a_{i+3}(i+1), \ldots, a_{n-1}(i+1)),
  \end{split}\]
}
and hence,
\begin{eqnarray}
  & & a_j(i) = a_j(i+1), \hspace{1em} 1 \leq j \leq i-1 \,\, \text{and} \,\, i+3 \leq j \leq n-1, \label{eq-tfc-1} \\
  & & a_{i+2}(i) = a_i(i+1), \label{eq-tfc-2} \\
  & & a_{i+1}(i) -a_{i+2}(i)+a_i(i) = a_{i+1}(i+1) -a_i(i+1) +a_{i+2}(i+1). \label{eq-tfc-3}
\end{eqnarray}
Similarly, for $i=n-2$, we obtain
\begin{eqnarray}
  & & a_j(n-2) = a_j(n-1), \hspace{1em} 1 \leq j \leq n-3, \label{eq-tfc-4} \\
  & & -|a(n-2)| = a_{n-2}(n-1), \label{eq-tfc-5} \\
  & & a_{n-1}(n-2) +|a(n-2)|+a_{n-2}(n-2) \nonumber \\
  & & \hspace{3em}  = a_{n-1}(n-1) -a_{n-2}(n-1) -|a(n-1)| \label{eq-tfc-6}
\end{eqnarray}
where $|a(i)|=a_1(i)+ \cdots + a_{n-1}(i)$ for any $1 \leq i \leq n-2$.
Similarly, we obtain the following equations.

\vspace{0.5em}

\noindent
{\bf{(B2)}}: We have
\[ (1- \sigma_j)f(\sigma_i)=(1-\sigma_i)f(\sigma_j). \]
For any $i<j \leq n-2$,
\begin{eqnarray}
  & & a_i(j) = a_{i+1}(j), \hspace{1em} i \leq j-2, \label{eq-tfc-7} \\
  & & a_j(i) = a_{j+1}(i), \hspace{1em} i \leq j-2. \label{eq-tfc-8} 
\end{eqnarray}
For any $i<j=n-1$,
\begin{eqnarray}
  & & |a(i)|+a_{n-1}(i) = 0, \hspace{1em} i \leq n-3, \label{eq-tfc-9} \\
  & & a_i(n-1) = a_{i+1}(n-1), \hspace{1em} i \leq n-3. \label{eq-tfc-10} 
\end{eqnarray}

\vspace{0.5em}

Here we sum up the above equations.
From (\ref{eq-tfc-1}) and (\ref{eq-tfc-4}),
\[ a_i(i+1)=a_i(i+2)= \cdots = a_i(n-1) \]
for any $1 \leq i \leq n-2$, and from (\ref{eq-tfc-1})
\[ a_i(1)=a_i(2)= \cdots =a_i(i-2) \]
for any $3 \leq i \leq n-1$. Moreover, from (\ref{eq-tfc-2}),
\[ a_i(i+1)=a_{i+2}(i). \]
This shows that all $a_i(j)$s except for $j =i, i-1$ are equal. Denote it by $\alpha$.
From (\ref{eq-tfc-9}), $|a(i)|=-\alpha$ for $1 \leq i \leq n-3$. From (\ref{eq-tfc-3}), we see
\[ a_i(i)+a_{i+1}(i) = -(n-2) \alpha \]
for any $1 \leq i \leq n-2$. Thus we see that $f(\sigma_i)$ for any $1 \leq i \leq n-1$ is determined by
\[ \alpha, a_1(1), a_2(2), \ldots, a_{n-1}(n-1), \]
and that there is no extra relations among them from the above equations.

\vspace{0.5em}

By the same argument as above, from {\bf{(SY2)}} and {\bf{(SY3)}},
if we set $\beta:=b_1(2)$ then we see that
$\beta = b_i(j)$ for $j \leq i-2$ or $i+1 \leq j$, $b_i(i)+b_{i+1}(i) = -(n-2) \beta$ for any $1 \leq i \leq n-2$,
and hence $f(s_i)$ for any $1 \leq i \leq n-1$ is determined by
\[ \beta, b_1(1), b_2(2), \ldots, b_{n-1}(n-1). \]
Consider {\bf{(SY1)}}. For $1 \leq i \leq n-2$, we obtain $2b_j(i)=0$ for $j \leq i-2$ and $b_i(i)+b_{i+1}(i)=0$ for $1 \leq i \leq n-2$.
For $i=n-1$, we obtain $b_j(n-1)=0$ for $1 \leq j \leq n-2$.
These conditions are equivalent to $\beta=0$.

\vspace{0.5em}

Consider {\bf{(BP1)}}. For any $i<j \leq n-2$,
\begin{eqnarray}
  & & b_i(j) = b_{i+1}(j), \hspace{1em} i \leq j-2, \nonumber \\
  & & a_j(i) = a_{j+1}(i), \hspace{1em} i \leq j-2. \nonumber 
\end{eqnarray}
For any $i<j=n-1$,
\begin{eqnarray}
  & & |a(i)|+a_{n-1}(i) = 0, \hspace{1em} i \leq n-3, \nonumber \\
  & & b_i(n-1) = b_{i+1}(n-1), \hspace{1em} i \leq n-3. \nonumber 
\end{eqnarray}
From these equations, we does not obtain a new condition among $a_i(j)$s and $b_i(j)$s.

\vspace{0.5em}

Consider {\bf{(BP2)}}. For any $i \leq n-3$,
\begin{eqnarray}
  & & a_{i+1}(i+1) - a_i(i)= b_{i+1}(i)-b_{i+2}(i+1), \label{eq-tfc-11} \\
  & & a_{i+2}(i+1) - a_{i+1}(i)= b_{i+2}(i+1)-b_{i+1}(i). \nonumber 
\end{eqnarray}
For the case of $i=n-2$,
\begin{equation}
  a_{n-1}(n-1) - a_{n-2}(n-2)= b_{n-1}(n-1)+b_{n-1}(n-2).  \label{eq-tfc-12}
\end{equation}
Thus we obtain a new condition from (\ref{eq-tfc-11}) and (\ref{eq-tfc-12}) which is equivalent to
\[ a_{i+1}(i+1) - a_i(i)= b_{i+1}(i+1)-b_i(i) \]
for $1 \leq i \leq n-2$.
From {\bf{(BP3)}}, we does not obtain a new condition among $a_i(j)$s and $b_i(j)$s.

\vspace{0.5em}

Therefore, we conclude that the homomorphism $\rho : Z^1(\mathrm{BP}_n,V) \rightarrow \Z^{\oplus n+1}$ defined by
$f \mapsto (\alpha, a_1(1), \ldots, a_{n-1}(n-1), b_1(1))$ is an isomorphism.
By combining the result for $B^1(\mathrm{BP}_n,V)$, in order to compute $H^1(\mathrm{BP}_n,V)$ through $\rho$,
it sufices to calculate the elementary divisors of the matrix
{\small
\[\begin{split} \bordermatrix{
               &  \alpha  & a_1(1) & a_2(2) & \cdots & \cdots & \cdots & a_{n-2}(n-2)  & a_{n-1}(n-1) & b_1(1) \\
           a_1 &  0       & -1     & 0      & \cdots & \cdots & \cdots & 0             & -1           & -1 \\
           a_2 &  \vdots  & 1      & -1     & 0      &        &        & \vdots        & -1           & 1 \\
           a_3 &  \vdots  & 0      & 1      & \ddots & \ddots &        & \vdots        & \vdots       & 0 \\
           a_4 &  \vdots  & \vdots & 0      & \ddots & \ddots & \ddots & \vdots        & \vdots       & \vdots \\
\,\,\, \vdots  &  \vdots  & \vdots & \vdots & 0      & \ddots & \ddots & 0             & \vdots       & \vdots \\
       a_{n-2} &  \vdots  & \vdots & \vdots & \vdots & \ddots & \ddots & -1            & -1           & \vdots \\
       a_{n-1} &  0       & 0      & 0      & 0      & \cdots & 0      & 1             & -2           & 0 \\
               }.  \end{split}\]
}
From this, we obtain the required result.
$\square$

\vspace{0.5em}

We remark that if we consider the braid group $B_n$ and the symmetric group $\mathfrak{S}_n$ as subgroups of $\mathrm{BP}_n$,
then by the same argument as the proof of Theorem {\rmfamily \ref{T-main-6}} we obtain the following results.

\begin{cor}\label{C-0617}
For $n \geq 3$, we have
\[ H^1(B_n,V) = \Z \oplus \Z/4\Z, \hspace{1em} H^1(\mathfrak{S}_n,V) = \Z/4\Z. \]
\end{cor}

\section{Acknowledgments}
The authors wolud like to express their sincere gratitude to Professor Yusuke Kuno for valuable comments for Lie algebras of tangential and
special derivations of the free Lie algebras, and Professor Toshiyuki Akita for valuable comments for the cohomology rings of the McCool groups.
We also would like to express their sincere gratitude to Professor Jaques Darn\'{e} for valuable comments for the McCool groups
in private communications.

\vspace{0.5em}

The first author is supported by JSPS KAKENHI Grant Number 18K03204.
The second author is supported by JSPS KAKENHI Grant Number 22K03299.

\end{document}